\documentclass{amsart}
\usepackage{graphicx}
\usepackage[all]{xy}
\usepackage{latexsym}
\usepackage{amssymb}
\usepackage{amsmath}

\xyoption{arc}

 \newfam\cyrfam

  \font\tencyr=wncyr10

  \font\sevencyr=wncyr7

  \font\fivecyr=wncyr5

  \textfont\cyrfam=\tencyr \scriptfont\cyrfam=\sevencyr

    \scriptscriptfont\cyrfam=\fivecyr

  \newfam\cyifam

  \font\tencyi=wncyi10

  \font\sevencyi=wncyi7

  \font\fivecyi=wncyi5

  \textfont\cyifam=\tencyi \scriptfont\cyifam=\sevencyi

    \scriptscriptfont\cyifam=\fivecyi

\def\id{{\mbox{1 \hskip -7pt 1}}}

 \newcommand{\lon}{\longrightarrow}
 \newcommand{\bu}{\bullet}
 
 \newcommand{\rar}{\rightarrow}
 \newcommand{\hook}{\hookrightarrow}

\newcommand{\p}{{\partial}}

\newcommand{\bubu}{{_{\xy
 (0,0)*{_\bullet},
(2.5,0)*{_\bu},
 (-0.5,-0.05)*{}="a",
(2.8,-0.05)*{}="b",
\ar @{-} "a";"b" <0pt>
\endxy}}}
\newcommand{\wibu}{{_{\xy
 (0,0)*{_\circ},
(2.5,0)*{_\bu},
 (0.26,-0.05)*{}="a",
(2.8,-0.05)*{}="b",
\ar @{-} "a";"b" <0pt>
\endxy}}}
\newcommand{\wiwi}{{_{\xy
 (0,0)*{_\circ},
(1.5,0)*{_\circ},
\endxy}}}

 \newcommand{\Z}{{\mathbb Z}}
 \newcommand{\bS}{{\mathbb S}}
 
 \newcommand{\C}{{\mathbb C}}
 \newcommand{\R}{{\mathbb R}}
 
 \newcommand{\K}{{\mathbb K}}

 \newcommand{\bbH}{{\mathbb H}}
\newcommand{\Conf}{{\mathit{Conf}}}
 \newcommand{\ot}{\otimes}

\newcommand{\sC}{{\mathsf C}}

\newcommand{\sG}{{\mathsf G}}
\newcommand{\sM}{{\mathsf M}}

\newcommand{\sa}{{\mathsf a}}
\newcommand{\ssf}{{\mathsf f}}

\newcommand{\Def}{\mathsf{Def}}

 \newcommand{\CoDer}{\mathsf{CoDer}}

 \newcommand{\Beq}{\begin{equation}}
 \newcommand{\Eeq}{\end{equation}}
 \newcommand{\Beqr}{\begin{eqnarray}}
 \newcommand{\Eeqr}{\end{eqnarray}}
 \newcommand{\Beqrn}{\begin{eqnarray*}}
 \newcommand{\Eeqrn}{\end{eqnarray*}}
 \newcommand{\Ba}{\begin{array}}
 \newcommand{\Ea}{\end{array}}
 \newcommand{\Bi}{\begin{itemize}}
 \newcommand{\Ei}{\end{itemize}}
 \newcommand{\Bc}{\begin{center}}
 \newcommand{\Ec}{\end{center}}
 \newcommand{\fg}{{\mathfrak g}}

\newcommand{\fr}{{\mathfrak r}}
\newcommand{\fs}{{\mathfrak s}}

\newcommand{\ft}{{\mathfrak t}}

\newcommand{\ii}{{\mathfrak i}}

\newcommand{\fW}{{\mathfrak W}}

 \newcommand{\f}{{\mathcal O}}
 \newcommand{\cA}{{\mathcal A}}
 
 \newcommand{\cC}{{\mathcal C}}
 
 \newcommand{\cE}{{\mathcal E}}
 \newcommand{\cF}{{\mathcal F}}
 \newcommand{\cG}{{\mathcal G}}

 \newcommand{\caL}{{\mathcal L}}
 \newcommand{\cM}{{\mathcal M}}
 
 \newcommand{\cP}{{\mathcal P}}

 \newcommand{\cT}{{\mathcal T}}

 \newcommand{\al}{\alpha}
 \newcommand{\be}{\beta}
 \newcommand{\ga}{\gamma}
 
 \newcommand{\Ga}{\Gamma}

 \newcommand{\var}{\varepsilon}
 
 \newcommand{\om}{\omega}

 \newcommand{\Ker}{{\mathsf K \mathsf e \mathsf r}\, }
 
 \newcommand{\Hom}{{\mathrm H\mathrm o\mathrm m}}
 
 \newcommand{\sip}{\smallskip}
 \newcommand{\bip}{\bigskip}
 \newcommand{\mip}{\vspace{2.5mm}}

\theoremstyle{plain}
\swapnumbers
\newtheorem{theorem}{Theorem}[subsection]
\newtheorem{corollary}[theorem]{Corollary}

\newtheorem{lemma}[theorem]{Lemma}
\newtheorem{proposition}[theorem]{Proposition}

\newtheorem{prop-def}[theorem]{Proposition-definition}

\newtheorem{f-theorem}{Formality Theorem}[section]
\newtheorem{main-theorem}{Main~Theorem}[section]
\newtheorem{section-theorem}{Theorem}[section]

\theoremstyle{definition}

\begin{document}

 \sloppy

 \newenvironment{proo}{\begin{trivlist} \item{\sc {Proof.}}}
  {\hfill $\square$ \end{trivlist}}

\long\def\symbolfootnote[#1]#2{\begingroup%
\def\thefootnote{\fnsymbol{footnote}}\footnote[#1]{#2}\endgroup}


\title{Grothendieck-Teichm\"uller group and
 Poisson cohomologies}

 \author{ Johan Alm and Sergei Merkulov}
\address{ Department of Mathematics, Stockholm University, 10691 Stockholm, Sweden}
\email{alm@math.su.se, sm@math.su.se}

 \maketitle
\markboth{J.\ Alm and S.A.\ Merkulov}{$GRT$ and Poisson cohohomology}

{\large
\section{\bf Introduction}
}

It is proven in \cite{Ta1,Wi} that the Grothendieck-Teichm\"uller group, $GRT$, acts up to homotopy
on the set, $\{\pi\}$, of Poisson structures (depending on a formal parameter $\hbar$) on an arbitrary smooth manifold.
Universal formulae for such an action can be represented as sums over Feynman graphs with
 weights given by integrals  over compactified configuration spaces introduced in \cite{Me-Auto,Me4}.

\sip

 Any Poisson structure makes the algebra of polyvector fields, $(\cT_{poly}(M)[[\hbar]], \wedge)$
into a {\em Poisson complex}, more precisely, into a differential graded (dg, for short) associative algebra
with the differential $d_\pi=[\pi,\ ]_S$, where $[\ ,\ ]_S$ is the Schouten bracket. The cohomology of this complex is
sometimes
denoted by $H^\bu(M, \pi)$ and is called the {\em Poisson cohomology}\ of $(M,\pi)$. The dg algebra
$(\cT_{poly}(M)[[\hbar]], \wedge, d_\pi)$ is of special type --- both operations $\wedge$ and $d_\pi$
respect the Schouten bracket in the sense of {\em dg Gerstenhaber algebra}. The main purpose of our paper is to study
\Bi
\item a class of universal $\cA ss_\infty$ structures on $\cT_{poly}(\R^d)$ which are consistent with
the Schouten bracket in the sense of   strong homotopy {\em non-commutative}\,  Gerstenhaber ($nc\cG_\infty$, for short) algebras;
\item universal actions of the group $GRT$ on this class,
\Ei
and then use these technical gadgets to give a constructive proof of the following
\mip

{\bf Main Theorem}. {\em  Let $\pi$ be a Poisson structure on $M$, $\ga$ an
arbitrary element
of $GRT$, and  let $\ga(\pi)$ be the Poisson structure on $M$ obtained from $\pi$ by an action of $\ga$. Then  there exists a morphism,
\Beq\label{1: F-gamma}
F^\ga: H^\bu( M, \pi) \lon H^\bu (M, \ga(\pi)),
\Eeq
 of  associative algebras.
}

\bip

The morphism (\ref{1: F-gamma}) is, in general, highly non-trivial.  In one of the simplest cases, when
$\pi$ is a linear Poisson structure on an affine manifold $M=\R^d$ (which is equivalent to the structure of a Lie algebra on the dual vector space $\fg:= (\R^d)^*$), the morphism $F^\ga$ becomes an
algebra {\em automorphism}\, of the Chevalley-Eilenberg cohomology of the $\fg$-module $\odot^\bu \fg:= \oplus_{n\geq 0} \odot^n \fg$,
$$
H^\bu(\R^d, \pi) = H^\bu( \R^d, \ga(\pi)) = H^\bu(\fg, \odot^\bu \fg),
$$
and its restriction to $H^0(\fg, \odot^\bu \fg)=(\odot^\bu \fg)^{\fg}$ coincides precisely with
 Kontsevich's  generalization of the classical Duflo map (see Theorems 7 and  8 in  \S 4.8 of \cite{Ko3}).
Thus in this special case our main theorem extends Kontsevich's action of $GRT$ on $(\odot^\bu \fg)^{\fg}$ to the full cohomology  $H^\bu(\fg, \odot^\bu \fg)$,
 and also gives explicit formulae for that extension.

 \mip

 The existence of the {\em algebra}\,  morphism (\ref{1: F-gamma}) is far from obvious. We prove in this paper a kind of
 ``no-go" theorem which says that there does {\em  not}\, exist a
{\em universal}\, (i.e.\ given by formulae applicable to any Poisson structure)  $\cA ss_\infty$-morphism of dg associative algebras,
$$
 \left(\cT_{poly}(M)[[\hbar]], \wedge, d_\pi\right) \lon
 \left(\cT_{poly}(M)[[\hbar]], \wedge, d_{\ga(\pi)}\right).
$$
Hence the algebra morphism (\ref{1: F-gamma}) can not be lifted to the level
of the associated Poisson complexes in such a way that the wedge multiplication is respected in the strong homotopy sense.

\mip

Our main technical tool is the deformation theory of universal $nc\cG_\infty$-structures on polyvector fields. This is governed by the mapping cone of a natural morphism
of graph complexes introduced and studied by Thomas Willwacher in \cite{Wi}. Using some of his results we show that there exists an exotic universal $GRT$-deformation of the standard the dg associative algebra $(\cT_{poly}(M), \wedge, d_{\ga(\pi)})$ into an $\cA ss_\infty$-algebra,
$$
\left(\cT_{poly}(M)[[\hbar]],\
 \mu_\bu^{\ga}=\{\mu_n^{\ga}\}_{n\geq 1}\right),\ \ \ga\in GRT,
$$
whose differential $\mu_1^\ga$ is independent of $\ga$ and equals  $d_{\pi}$, while the higher homotopy
operations $\mu^\ga_{n\geq 2}$ are independent of $\pi$ and are fully determined by $\wedge$ and $\ga$. This universal $\cA ss_\infty$- algebra structure on $\cT_{poly}(M)$ is homotopy equivalent to  $(\cT_{poly}(M), \wedge, d_{\ga(\pi)})$, i.e.\ there exists
a universal continuous   $\cA ss_\infty$ isomorphism,
\Beq\label{1: cF-gamma}
\cF^\ga:  \left(\cT_{poly}(M)[[\hbar]], \mu_\bu^{\ga}\right) \lon  \left(\cT_{poly}(M)[[\hbar]], \wedge, d_{\ga(\pi)}\right)
\Eeq
which on cohomology induces the map (\ref{1: F-gamma}) and hence proves the main theorem.

\mip

The above mentioned deformation,  $\left(\cT_{poly}(M)[[\hbar]], \mu_\bu^{\ga}\right)$, of the standard dg algebra structure on the space of polyvector fields
on a Poisson manifold is of a special type --- it respects the Schouten bracket.  Omitting reference to a
particular Poisson structure  $\pi$ on $M$, we can say that we study in this paper universal deformations
of the standard Gerstenhaber algebra structure on polyvector fields in the class
of $nc\cG_\infty$-algebras (rather than in the class of $\cG_\infty$-algebras). The 2-coloured operad $nc\cG_\infty$ is a minimal resolution of the 2-coloured  Koszul operad, $nc\cG$, of non-commutative Gerstenhaber algebras, and, moreover, it admits a very natural geometric realization via configuration spaces of points in the pair, $\R\subset \C$, consisting of the complex plane $\C$ and a line
$\R$ drawn in the plane \cite{A}.
 We prove that, up to $nc\cG_\infty$-isomorphisms there are only two universal $nc\cG_\infty$-structures on polyvector fields, the one which comes from the standard Gerstenhaber algebra structure, and the exotic one which was
introduced in \cite{A} in terms of a de Rham field theory on a certain operad of compactified
configuration spaces.



\sip

\subsection{Some notation}
 The set $\{1,2, \ldots, n\}$ is abbreviated to $[n]$;  its group of automorphisms is
denoted by $\bS_n$. The
cardinality of a finite set
$A$ is denoted by $\# A$. If $V=\oplus_{i\in \Z} V^i$ is a graded vector space, then
$V[k]$ stands for the graded vector space with $V[k]^i:=V^{i+k}$ and
and $s^k$ for the associated isomorphism $V\rar V[k]$; for $v\in V^i$ we set $|v|:=i$.
For a pair of graded vector spaces $V_1$ and $V_2$, the symbol $\Hom_i(V_1,V_2)$ stands for the space
of homogeneous linear maps of degree $i$, and $\Hom(V_1,V_2):=\bigoplus_{i\in \Z}\Hom_i(V_1,V_2)$; for example,
$s^k\in \Hom_{-k}(V,V[k])$. For an
operad $\cP$ we denote by $\cP\{k\}$ the unique operad which has the following property:
for any graded vector space $V$ there is a one-to-one correspondence between representations of
$\cP\{k\}$ in $V$ and representations of
$\cP$ in $V[-k]$; in particular, $\cE nv_V\{k\}=\cE nd_{V[k]}$.

\subsection{Acknowledgements} {It is a great pleasure to thank
Thomas Willwacher  for many very useful discussions and correspondences.}

\bip

{\large
\section{\bf Compactified configuration spaces of points in the flag $\R\subset \C$}
}

\subsection{$ \caL ie_\infty$-algebras} For a finite set $A$ let $\Conf_A(\C)$ stand for the set
of all injections, $\{A\hook \C\}$. For $\# A \geq 2$ the orbit space
$$
C_A(\C):= \frac{\{A\hook \C\}}{z\rar \R^+ z + \C},
$$
is naturally a real $(2\# A-3)$-dimensional manifold (if $A=[n]$, we use the notations $C_n(\C)$).
Its Fulton-MacPherson compactification, $\overline{C}_A(\C)$, can be made into a compact smooth manifold with
corners \cite{Ko2} (or into
a compact semialgebraic manifold). Moreover, the collection
$$
\overline{C}(\C)=\{\overline{C}_A(\C)\}_{\# A \geq 2}
$$
has a natural structure of a non-unital pseudo-operad in the category of oriented smooth manifolds with corners.
The associated operad of chains, $\cC hains( \overline{C}(\C))$, contains a suboperad of fundamental chains,
 $\cF \cC hains( \overline{C}(\C))$, which is precisely the operad, $\caL_\infty\{1\}$, of degree shifted $L_\infty$-algebras
 (see \cite{Me3} for a review).

\subsection{OCHA versus strong homotopy non-commutative Gerstenhaber algebras} For arbitrary finite sets $A$ and $B$ consider the space of injections,
$$
\Conf_{A,B}(\C):=\{  A\sqcup B \hook \C, B\hook \R\subset \C\},
$$
and, for $2\# A + \# B\geq 2$, consider the quotient space,
$$
C_{A,B}(\C):= \frac{\Conf_{A,B}(\C)}{z\rar \R^+ z + \R},
$$
by the affine group $\R^+ \ltimes \R$. As $\C\setminus \R= \bbH \sqcup \bbH^-$, where
$\bbH $ (resp. $\bbH^-$) is the upper (resp., lower) upper half-plane, we can consider subspaces,
$$
\Conf_{A,B}(\bbH):=\{  A\hook \bbH, B\hook \R\} \subset \Conf_{A,B}(\C)
$$
and
$$
C_{A,B}(\bbH):= \frac{\Conf_{A,B}(\bbH)}{z\rar \R^+ z + \R}  \subset C_{A,B}(\C)
$$
\mip

\[
\resizebox{5cm}{!}
{
\xy
(-16,0)*{\circ},
(-16,2)*{^{x_1}},
(20,0)*{\circ},
(20,2)*{^{x_2}},
(-8,23)*{\bullet}="e",
(7.2,-10)*{\bullet}="f",
(-10,24.5)*{_{z_1}},
(9.5,-8.5)*{_{z_2}},
(-10,-11)*{\bullet},
(-11.8,-9.8)*{^{z_3}},
(-26,16)*{\bullet},
(-27.7,16.8)*{^{z_5}},
(26,16)*{\bullet},
(27.7,16.8)*{^{z_4}},
(-2.5,6)*{\cdot},
(-2.5,6)*{\ \ \ \ ^{\ii}},
 (-35,0)*{}="a",
(35,0)*{}="b",
(-2.5,-15)*{}="c",
(-2.5,30)*{}="d",
\ar @{->} "a";"b" <0pt>
\ar @{.>} "c";"d" <0pt>
\endxy}
\hspace{33mm}
\resizebox{5.2cm}{!}{
\xy
(-16,0)*{\circ},
(-16,2)*{^{x_1}},
(20,0)*{\circ},
(20,2)*{^{x_2}},
(-8,23)*{\bullet}="e",
(7.2,10)*{\bullet}="f",
(-10,24.5)*{_{z_1}},
(9.5,8.5)*{_{z_2}},
(-10,11)*{\bullet},
(-9.7,12.8)*{^{z_3}},
(-26,16)*{\bullet},
(-27.7,16.8)*{^{z_n}},
(26,16)*{\bullet},
(27.7,16.8)*{^{z_4}},
(-2.5,6)*{\cdot},
(-2.5,6)*{\ \ \ \ ^{\ii}},
 (-40,0)*{}="a",
(40,0)*{}="b",
(-2.5,0)*{}="c",
(-2.5,30)*{}="d",
\ar @{->} "a";"b" <0pt>
\ar @{.>} "c";"d" <0pt>
\endxy}
\]
$$
C_{5,2}(\C) \hspace{70mm} C_{5,2}(\bbH) \ \ \
$$
\mip

The Fulton-MacPherson compactification, $\overline{C}_{A,B}(\bbH)$, of $C_{A,B}(\bbH)$ was introduced in \cite{Ko2}.
The fundamental chain complex, $\cF\cC hains(\overline{C}(\bbH))$, of the disjoint union,
$$
\overline{C}(\bbH):=\overline{C}_\bu(\C)\bigsqcup \overline{C}_{\bu,\bu}(\bbH),
$$
is a dg quasi-free 2-colored operad \cite{KS} generated by
 \Bi
 \item[(i)] degree $3-2n$ corollas,
\Beq\label{2: Lie_inf corolla}
\resizebox{2.2cm}{!}{
\xy
(1,-5)*{\ldots},
(-13,-7)*{_1},
(-8,-7)*{_2},
(-3,-7)*{_3},
(7,-7)*{_{n-1}},
(13,-7)*{_n},
 (0,0)*{\bu}="a",
(0,5)*{}="0",
(-12,-5)*{}="b_1",
(-8,-5)*{}="b_2",
(-3,-5)*{}="b_3",
(8,-5)*{}="b_4",
(12,-5)*{}="b_5",
\ar @{-} "a";"0" <0pt>
\ar @{-} "a";"b_2" <0pt>
\ar @{-} "a";"b_3" <0pt>
\ar @{-} "a";"b_1" <0pt>
\ar @{-} "a";"b_4" <0pt>
\ar @{-} "a";"b_5" <0pt>
\endxy}
=
\resizebox{2.2cm}{!}{
\xy
(1,-6)*{\ldots},
(-13,-7)*{_{\sigma(1)}},
(-6.7,-7)*{_{\sigma(2)}},
(13,-7)*{_{\sigma(n)}},
 (0,0)*{\bu}="a",
(0,5)*{}="0",
(-12,-5)*{}="b_1",
(-8,-5)*{}="b_2",
(-3,-5)*{}="b_3",
(8,-5)*{}="b_4",
(12,-5)*{}="b_5",
\ar @{-} "a";"0" <0pt>
\ar @{-} "a";"b_2" <0pt>
\ar @{-} "a";"b_3" <0pt>
\ar @{-} "a";"b_1" <0pt>
\ar @{-} "a";"b_4" <0pt>
\ar @{-} "a";"b_5" <0pt>
\endxy},
\ \ \ \forall \sigma\in \bS_n,\ n\geq2
\Eeq
representing $\overline{C}_n(\C)$,  and
\item[(ii)]
degree $2-2n-m$ corollas,
\Beq\label{2: OCHA corolla}
\Ba{c}
\resizebox{2.7cm}{!}{
\begin{xy}
 <0mm,-0.5mm>*{\blacktriangledown};
 <0mm,0mm>*{};<0mm,5mm>*{}**@{.},
 <0mm,0mm>*{};<-16mm,-5mm>*{}**@{-},
 <0mm,0mm>*{};<-11mm,-5mm>*{}**@{-},
 <0mm,0mm>*{};<-3.5mm,-5mm>*{}**@{-},
 <0mm,0mm>*{};<-6mm,-5mm>*{...}**@{},
   <0mm,0mm>*{};<-16mm,-8mm>*{^{1}}**@{},
   <0mm,0mm>*{};<-11mm,-8mm>*{^{2}}**@{},
   <0mm,0mm>*{};<-3mm,-8mm>*{^{n}}**@{},
 <0mm,0mm>*{};<16mm,-5mm>*{}**@{.},
 <0mm,0mm>*{};<8mm,-5mm>*{}**@{.},
 <0mm,0mm>*{};<3.5mm,-5mm>*{}**@{.},
 <0mm,0mm>*{};<11.6mm,-5mm>*{...}**@{},
   <0mm,0mm>*{};<19mm,-8mm>*{^{\bar{m}}}**@{},
<0mm,0mm>*{};<10mm,-8mm>*{^{\bar{2}}}**@{},
   <0mm,0mm>*{};<5mm,-8mm>*{^{\bar{1}}}**@{},
 \end{xy}}
\Ea=
\Ba{c}
\resizebox{2.7cm}{!}{
\begin{xy}
 <0mm,-0.5mm>*{\blacktriangledown};
 <0mm,0mm>*{};<0mm,5mm>*{}**@{.},
 <0mm,0mm>*{};<-16mm,-5mm>*{}**@{-},
 <0mm,0mm>*{};<-11mm,-5mm>*{}**@{-},
 <0mm,0mm>*{};<-3.5mm,-5mm>*{}**@{-},
 <0mm,0mm>*{};<-6mm,-5mm>*{...}**@{},
   <0mm,0mm>*{};<-18mm,-8mm>*{^{\sigma(1)}}**@{},
   <0mm,0mm>*{};<-11mm,-8mm>*{^{\sigma(2)}}**@{},
   <0mm,0mm>*{};<-3mm,-8mm>*{^{\sigma(n)}}**@{},
 <0mm,0mm>*{};<16mm,-5mm>*{}**@{.},
 <0mm,0mm>*{};<8mm,-5mm>*{}**@{.},
 <0mm,0mm>*{};<3.5mm,-5mm>*{}**@{.},
 <0mm,0mm>*{};<11.6mm,-5mm>*{...}**@{},
   <0mm,0mm>*{};<19mm,-8mm>*{^{\bar{m}}}**@{},
<0mm,0mm>*{};<10mm,-8mm>*{^{\bar{2}}}**@{},
   <0mm,0mm>*{};<5mm,-8mm>*{^{\bar{1}}}**@{},
 \end{xy}}
\Ea
, \ \ \ \ 2n+m\geq 2, \forall\ \sigma\in \bS_n
\Eeq
representing  $\overline{C}_{n,m}(\bbH)$.
\Ei
The differential  in  $\cF\cC hains(\overline{C}(\bbH))$ is given on the generators
by \cite{Ko2, KS}
\Beqr
\p
\Ba{c}
\resizebox{2.5cm}{!}{
 \xy
(1,-5)*{\ldots},
(-13,-7)*{_1},
(-8,-7)*{_2},
(-3,-7)*{_3},
(7,-7)*{_{n-1}},
(13,-7)*{_n},
 (0,0)*{\bu}="a",
(0,5)*{}="0",
(-12,-5)*{}="b_1",
(-8,-5)*{}="b_2",
(-3,-5)*{}="b_3",
(8,-5)*{}="b_4",
(12,-5)*{}="b_5",
\ar @{-} "a";"0" <0pt>
\ar @{-} "a";"b_2" <0pt>
\ar @{-} "a";"b_3" <0pt>
\ar @{-} "a";"b_1" <0pt>
\ar @{-} "a";"b_4" <0pt>
\ar @{-} "a";"b_5" <0pt>
\endxy}
\Ea
&=&
\sum_{A\varsubsetneq [n]\atop
\# A\geq 2}\hspace{-2mm}
\Ba{c}
\resizebox{2.4cm}{!}{
\begin{xy}
<10mm,0mm>*{\bu},
<10mm,0.8mm>*{};<10mm,5mm>*{}**@{-},
<0mm,-10mm>*{...},
<14mm,-5mm>*{\ldots},
<13mm,-7mm>*{\underbrace{\ \ \ \ \ \ \ \ \ \ \ \ \  }},
<14mm,-10mm>*{_{[n]\setminus A}};
<10.3mm,0.1mm>*{};<20mm,-5mm>*{}**@{-},
<9.7mm,-0.5mm>*{};<6mm,-5mm>*{}**@{-},
<9.9mm,-0.5mm>*{};<10mm,-5mm>*{}**@{-},
<9.6mm,0.1mm>*{};<0mm,-4.4mm>*{}**@{-},
<0mm,-5mm>*{\bu};
<-5mm,-10mm>*{}**@{-},
<-2.7mm,-10mm>*{}**@{-},
<2.7mm,-10mm>*{}**@{-},
<5mm,-10mm>*{}**@{-},
<0mm,-12mm>*{\underbrace{\ \ \ \ \ \ \ \ \ \ }},
<0mm,-15mm>*{_{A}},
\end{xy}}
\Ea
\label{2: d on Lie corollas}
\\
\p
\Ba{c}
\resizebox{2.6cm}{!}{
\begin{xy}
 <0mm,-0.5mm>*{\blacktriangledown};
 <0mm,0mm>*{};<0mm,5mm>*{}**@{.},
 <0mm,0mm>*{};<-16mm,-5mm>*{}**@{-},
 <0mm,0mm>*{};<-11mm,-5mm>*{}**@{-},
 <0mm,0mm>*{};<-3.5mm,-5mm>*{}**@{-},
 <0mm,0mm>*{};<-6mm,-5mm>*{...}**@{},
   <0mm,0mm>*{};<-16mm,-8mm>*{^{1}}**@{},
   <0mm,0mm>*{};<-11mm,-8mm>*{^{2}}**@{},
   <0mm,0mm>*{};<-3mm,-8mm>*{^{n}}**@{},
 <0mm,0mm>*{};<16mm,-5mm>*{}**@{.},
 <0mm,0mm>*{};<8mm,-5mm>*{}**@{.},
 <0mm,0mm>*{};<3.5mm,-5mm>*{}**@{.},
 <0mm,0mm>*{};<11.6mm,-5mm>*{...}**@{},
   <0mm,0mm>*{};<17mm,-8mm>*{^{\bar{m}}}**@{},
<0mm,0mm>*{};<10mm,-8mm>*{^{\bar{2}}}**@{},
   <0mm,0mm>*{};<5mm,-8mm>*{^{\bar{1}}}**@{},
 \end{xy}}
\Ea
&=& -
\sum_{A\varsubsetneq [n]\atop
\# A\geq 2}\
\Ba{c}
\resizebox{2.9cm}{!}{
\begin{xy}
 <0mm,-0.5mm>*{\blacktriangledown};
 <0mm,0mm>*{};<0mm,5mm>*{}**@{.},
 <0mm,0mm>*{};<-16mm,-5mm>*{}**@{-},
 <0mm,0mm>*{};<-11mm,-5mm>*{}**@{-},
 <0mm,0mm>*{};<-3.5mm,-5mm>*{}**@{-},
 <0mm,0mm>*{};<-6mm,-5mm>*{...}**@{},
 <0mm,0mm>*{};<16mm,-5mm>*{}**@{.},
 <0mm,0mm>*{};<8mm,-5mm>*{}**@{.},
 <0mm,0mm>*{};<3.5mm,-5mm>*{}**@{.},
 <0mm,0mm>*{};<11.6mm,-5mm>*{...}**@{},
   <0mm,0mm>*{};<17mm,-8mm>*{^{\bar{m}}}**@{},
<0mm,0mm>*{};<10mm,-8mm>*{^{\bar{2}}}**@{},
   <0mm,0mm>*{};<5mm,-8mm>*{^{\bar{1}}}**@{},
<-17mm,-12mm>*{\underbrace{\ \ \ \ \ \ \ \ \ \   }},
<-17mm,-14.9mm>*{_A};
<-6mm,-7mm>*{\underbrace{\ \ \ \ \ \ \  }},
<-6mm,-10mm>*{_{[n]\setminus A}};
 (-16.5,-5.5)*{\bu}="a",
(-23,-10)*{}="b_1",
(-20,-10)*{}="b_2",
(-16,-10)*{...}="b_3",
(-12,-10)*{}="b_4",
\ar @{-} "a";"b_2" <0pt>
\ar @{-} "a";"b_1" <0pt>
\ar @{-} "a";"b_4" <0pt>
 \end{xy}}
\Ea
\label{2: differential on Konts corollas}\\
&+& \sum_{k, l, [n]=I_1\sqcup I_2\atop
{2\# I_1 + m \geq l+1 \atop
2\#I_2 + l\geq 2}}
(-1)^{k+l(n-k-l)}
\Ba{c}
\resizebox{3.5cm}{!}{
\begin{xy}
 <0mm,-0.5mm>*{\blacktriangledown};
 <0mm,0mm>*{};<0mm,6mm>*{}**@{.},
 <0mm,0mm>*{};<-16mm,-6mm>*{}**@{-},
 <0mm,0mm>*{};<-11mm,-6mm>*{}**@{-},
 <0mm,0mm>*{};<-3.5mm,-6mm>*{}**@{-},
 <0mm,0mm>*{};<-6mm,-6mm>*{...}**@{},
<0mm,0mm>*{};<2mm,-9mm>*{^{\bar{1}}}**@{},
<0mm,0mm>*{};<6mm,-9mm>*{^{\bar{k}}}**@{},
<0mm,0mm>*{};<19mm,-9mm>*{^{\overline{k+l+1}}}**@{},
<0mm,0mm>*{};<28mm,-9mm>*{^{\overline{m}}}**@{},
<0mm,0mm>*{};<13mm,-16.6mm>*{^{\overline{k+1}}}**@{},
<0mm,0mm>*{};<20mm,-16.6mm>*{^{\overline{k+l}}}**@{},
 <0mm,0mm>*{};<11mm,-6mm>*{}**@{.},
 <0mm,0mm>*{};<6mm,-6mm>*{}**@{.},
 <0mm,0mm>*{};<2mm,-6mm>*{}**@{.},
 <0mm,0mm>*{};<17mm,-6mm>*{}**@{.},
 <0mm,0mm>*{};<25mm,-6mm>*{}**@{.},
 <0mm,0mm>*{};<4mm,-6mm>*{...}**@{},
<0mm,0mm>*{};<20mm,-6mm>*{...}**@{},
<6.5mm,-16mm>*{\underbrace{\ \ \ \ \   }_{I_2}},
<-10mm,-9mm>*{\underbrace{\ \ \ \ \ \ \ \ \ \ \ \   }_{I_1}},
 %
 (11,-7)*{\blacktriangledown}="a",
(4,-13)*{}="b_1",
(9,-13)*{}="b_2",
(16,-13)*{...},
(7,-13)*{...},
(13,-13)*{}="b_3",
(19,-13)*{}="b_4",
\ar @{-} "a";"b_2" <0pt>
\ar @{.} "a";"b_3" <0pt>
\ar @{-} "a";"b_1" <0pt>
\ar @{.} "a";"b_4" <0pt> \nonumber
 \end{xy}}
\Ea
\Eeqr
Representations of  $(\cF\cC hains(\overline{C}(\bbH)), \p)$ in a pair of dg vector spaces
$(A,\fg)$  were called in \cite{KS} {\em open-closed homotopy algebras}\, or OCHAs for short.
Such a representation, $\rho$, is uniquely determined
by its values on the generators,
\Beqrn
\nu_n&:=& \rho  \left(\Ba{c} \resizebox{2.5cm}{!}{  \xy
(1,-5)*{\ldots},
(-13,-7)*{_1},
(-8,-7)*{_2},
(-3,-7)*{_3},
(7,-7)*{_{n-1}},
(13,-7)*{_n},
 (0,0)*{\bu}="a",
(0,5)*{}="0",
(-12,-5)*{}="b_1",
(-8,-5)*{}="b_2",
(-3,-5)*{}="b_3",
(8,-5)*{}="b_4",
(12,-5)*{}="b_5",
\ar @{-} "a";"0" <0pt>
\ar @{-} "a";"b_2" <0pt>
\ar @{-} "a";"b_3" <0pt>
\ar @{-} "a";"b_1" <0pt>
\ar @{-} "a";"b_4" <0pt>
\ar @{-} "a";"b_5" <0pt>
\endxy}\Ea\right) \in \Hom(\fg^{\odot n}, \fg)[3-2n], \ \ \ n\geq 2, \\
\mu_{n,m} & := & \rho  \left(
\Ba{c}
\resizebox{2.7cm}{!}{
\begin{xy}
 <0mm,-0.5mm>*{\blacktriangledown};
 <0mm,0mm>*{};<0mm,5mm>*{}**@{.},
 <0mm,0mm>*{};<-16mm,-5mm>*{}**@{-},
 <0mm,0mm>*{};<-11mm,-5mm>*{}**@{-},
 <0mm,0mm>*{};<-3.5mm,-5mm>*{}**@{-},
 <0mm,0mm>*{};<-6mm,-5mm>*{...}**@{},
   <0mm,0mm>*{};<-16mm,-8mm>*{^{1}}**@{},
   <0mm,0mm>*{};<-11mm,-8mm>*{^{2}}**@{},
   <0mm,0mm>*{};<-3mm,-8mm>*{^{n}}**@{},
 <0mm,0mm>*{};<16mm,-5mm>*{}**@{.},
 <0mm,0mm>*{};<8mm,-5mm>*{}**@{.},
 <0mm,0mm>*{};<3.5mm,-5mm>*{}**@{.},
 <0mm,0mm>*{};<11.6mm,-5mm>*{...}**@{},
   <0mm,0mm>*{};<17mm,-8mm>*{^{\bar{m}}}**@{},
<0mm,0mm>*{};<10mm,-8mm>*{^{\bar{2}}}**@{},
   <0mm,0mm>*{};<5mm,-8mm>*{^{\bar{1}}}**@{},
 \end{xy}}
\Ea
\right) \in \Hom(\fg^{\odot n}\otimes A^{\ot m}, A)[2-2n-m], \ \ 2n +m \geq 2,
\Eeqrn
which satisfy quadratic relations given by the above formulae for the differential $\p$ and give us, therefore,
the following list of algebraic structures in $(A,\fg)$:
\Bi
\item[(i)] an $\caL_\infty\{1\}$-algebra structure, $\nu_\bu=\{\nu_n: \odot^n \fg\rar \fg[3-2n]\}_{n\geq 1}$, in $\fg$;
\item[(ii)] an $\cA_\infty$-algebra structure, $\mu_\bu=\{\mu_{0,m}: \ot^m A\rar A[2-m]\}_{m\geq 1}$, in $A$; if
$[\ ,\ ]_G$ stands for the standard Gerstenhaber bracket on the Hochschild cochains
$C(A,A)=\prod_{n\geq 0}\Hom (A^{\ot n}, A)[1-n])$,then $\mu_\bu$ defines a differential on $C(A,A)$, $d_\mu:=[\mu_\bu,\ ]_G$;
\item[(iii)] an  $\caL_\infty$-morphism, $F$, from the $L_\infty$-algebra $(\fg, \nu)$ to
the dg Lie algebra $(C(A,A), [\ ,\ ]_G, d_\mu)$.
\Ei

\sip

If $\rho$ is an arbitrary representation of $(\cF\cC hains(\overline{C}(\bbH)), \p)$  and
 $\ga\in \fg$ is an arbitrary Maurer-Cartan
element\footnote{We tacitly assume here that the $L_\infty$-algebra
$(X_c,\nu_\bu)$ is appropriately filtered so that
the MC equation makes sense. In our applications below $\nu_{n\geq 3}=0$
 so that one has no problems with convergence of the infinite sum.},
$$
\sum_{n\geq 0}\frac{1}{n!} \nu_n(\ga^{\ot n})=0, \ \ \ |\ga|=2,
$$
of the associated $\caL_\infty$-algebra $(\fg,\nu_\bu)$, then the maps
$$
\Ba{rccc}
\mu_m: & \ot^m A & \lon & A[[\hbar]][2-m], \ \ m\geq 0, \\
       & x_1\ot\ldots x_m & \lon &  \sum_{n\geq 1}\frac{\hbar^n}{n!}
\mu_{n,m}(\ga^{\ot n}\ot x_1\ot\ldots x_m)
\Ea
$$
make the topological (with respect to the adic topology) vector space $A[[\hbar]]$ into a
topological, {\em non-flat} (in general)\,  $\cA_\infty$-algebra (here $\hbar$ is a formal parameter, and
$A[[\hbar]:=A\ot\K[[\hbar]]$). Non-flatness originates from the generators (\ref{2: OCHA corolla})
with $m=0$, $n\geq 1$, which correspond to the boundary strata in $\overline{C}(\bbH)$ that are given by groups of points in the upper half plane collapsing to a point on the real line. It is clear how to get rid of such strata --- one should allow configurations of points everywhere in $\C$, and hence consider the Fulton-MacPherson compactifications \cite{A} of the configuration spaces
$C_{A,B}(\C)$ rather than $C_{A,B}(\bbH)$. The disjoint union
$$
\overline{CF}(\C):= \overline{C}_\bu(\C)\bigsqcup \overline{C}_{\bu,\bu}(\C),
$$
has a natural structure of a dg quasi-free 2-colored operad in the category of compact manifolds with corners. This operad is free in the category of sets. The suboperad,
$$
nc\cG_\infty:=\cF\cC hains(\overline{CF}(\C)),
$$
 of the associated chain operad  $\cC hains(\overline{CF}(\C))$ generated by fundamental chains is free in the category of graded vector spaces and
  is canonically isomorphic as a dg operad to the quotient operad
$$
nc\cG_\infty:= \cF\cC hains(\overline{C}(\bbH))/ I,
$$
where $I$ is the (differential) ideal generated by corollas (\ref{2: OCHA corolla}) with $m=0$, $n\geq 1$. The notation $nc\cG_\infty$ stems from the fact \cite{A} that this operad is a minimal resolution of a 2-coloured quadratic operad
which governs the type of algebras introduced in \cite{FGV} under the name of {\em Leibniz pairs}. Let us compare this
 quadratic operad with the operad, $\cG$, of Gerstenhaber algebras. The latter is a 1-colored quadratic
 operad generated by commutative associative product in degree 0,
$\Ba{c}\begin{xy}
 <0mm,0.65mm>*{};<0mm,3.5mm>*{}**@{-},
 <0.5mm,-0.5mm>*{};<2.2mm,-2.2mm>*{}**@{-},
 <-0.48mm,-0.48mm>*{};<-2.2mm,-2.2mm>*{}**@{-},
 <0mm,0mm>*{\circ};<0mm,0mm>*{}**@{},
 <0.5mm,-0.5mm>*{};<2.8mm,-3.3mm>*{_2}**@{},
 <-0.48mm,-0.48mm>*{};<-2.7mm,-3.3mm>*{_1}**@{},
 \end{xy}\Ea
=
\Ba{c}\begin{xy}
 <0mm,0.65mm>*{};<0mm,3.5mm>*{}**@{-},
 <0.5mm,-0.5mm>*{};<2.2mm,-2.2mm>*{}**@{-},
 <-0.48mm,-0.48mm>*{};<-2.2mm,-2.2mm>*{}**@{-},
 <0mm,0mm>*{\circ};<0mm,0mm>*{}**@{},
 <0.5mm,-0.5mm>*{};<2.8mm,-3.3mm>*{_1}**@{},
 <-0.48mm,-0.48mm>*{};<-2.7mm,-3.3mm>*{_2}**@{},
 \end{xy}\Ea$
and Lie bracket of degree $-1$,
$\Ba{c}
\begin{xy}
 <0mm,0.55mm>*{};<0mm,3.5mm>*{}**@{-},
 <0.5mm,-0.5mm>*{};<2.2mm,-2.2mm>*{}**@{-},
 <-0.48mm,-0.48mm>*{};<-2.2mm,-2.2mm>*{}**@{-},
 <0mm,0mm>*{\bu};<0mm,0mm>*{}**@{},
 <0.5mm,-0.5mm>*{};<2.8mm,-3.3mm>*{_2}**@{},
 <-0.48mm,-0.48mm>*{};<-2.7mm,-3.3mm>*{_1}**@{},
 \end{xy}\Ea
=
\Ba{c}
\begin{xy}
 <0mm,0.55mm>*{};<0mm,3.5mm>*{}**@{-},
 <0.5mm,-0.5mm>*{};<2.2mm,-2.2mm>*{}**@{-},
 <-0.48mm,-0.48mm>*{};<-2.2mm,-2.2mm>*{}**@{-},
 <0mm,0mm>*{\bu};<0mm,0mm>*{}**@{},
 <0.5mm,-0.5mm>*{};<2.8mm,-3.3mm>*{_1}**@{},
 <-0.48mm,-0.48mm>*{};<-2.7mm,-3.3mm>*{_2}**@{},
 \end{xy}\Ea$,
satisfying the compatibility condition
\Beq\label{2: Gertenshaber compatibility}
\Ba{c}
\begin{xy}
 <0mm,0mm>*{\circ};<0mm,0mm>*{}**@{},
 <0mm,0.49mm>*{};<0mm,3.5mm>*{}**@{-},
 <0.49mm,-0.49mm>*{};<1.9mm,-1.9mm>*{}**@{-},
 <-0.5mm,-0.5mm>*{};<-1.9mm,-1.9mm>*{}**@{-},
 <2.3mm,-2.3mm>*{\bu};<-2.3mm,-2.3mm>*{}**@{},
 <1.8mm,-2.8mm>*{};<0mm,-4.9mm>*{}**@{-},
 <2.8mm,-2.9mm>*{};<4.6mm,-4.9mm>*{}**@{-},
   <0.49mm,-0.49mm>*{};<-2.8mm,-3.2mm>*{_1}**@{},
   <-1.8mm,-2.8mm>*{};<0.1mm,-6.2mm>*{_2}**@{},
   <-2.8mm,-2.9mm>*{};<5.3mm,-6.2mm>*{_3}**@{},
 \end{xy}
 \Ea
\ = \
\Ba{c}
\begin{xy}
 <0mm,0mm>*{\circ};<0mm,0mm>*{}**@{},
 <0mm,0.49mm>*{};<0mm,3.5mm>*{}**@{-},
 <0.49mm,-0.49mm>*{};<1.9mm,-1.9mm>*{}**@{-},
 <-0.5mm,-0.5mm>*{};<-1.9mm,-1.9mm>*{}**@{-},
 <-2.3mm,-2.3mm>*{\bu};<-2.3mm,-2.3mm>*{}**@{},
 <-1.8mm,-2.8mm>*{};<0mm,-4.9mm>*{}**@{-},
 <-2.8mm,-2.9mm>*{};<-4.6mm,-4.9mm>*{}**@{-},
   <0.49mm,-0.49mm>*{};<2.7mm,-3.2mm>*{^3}**@{},
   <-1.8mm,-2.8mm>*{};<0.4mm,-6.2mm>*{^2}**@{},
   <-2.8mm,-2.9mm>*{};<-5.1mm,-6.2mm>*{^1}**@{},
 \end{xy}\Ea
\ + \
\Ba{c}
\begin{xy}
 <0mm,0mm>*{\circ};<0mm,0mm>*{}**@{},
 <0mm,0.49mm>*{};<0mm,3.5mm>*{}**@{-},
 <0.49mm,-0.49mm>*{};<1.9mm,-1.9mm>*{}**@{-},
 <-0.5mm,-0.5mm>*{};<-1.9mm,-1.9mm>*{}**@{-},
 <2.3mm,-2.3mm>*{\bu};<-2.3mm,-2.3mm>*{}**@{},
 <1.8mm,-2.8mm>*{};<0mm,-4.9mm>*{}**@{-},
 <2.8mm,-2.9mm>*{};<4.6mm,-4.9mm>*{}**@{-},
   <0.49mm,-0.49mm>*{};<-2.8mm,-3.2mm>*{_2}**@{},
   <-1.8mm,-2.8mm>*{};<0.1mm,-6.2mm>*{_1}**@{},
   <-2.8mm,-2.9mm>*{};<5.3mm,-6.2mm>*{_3}**@{},
 \end{xy}
 \Ea
\Eeq
This condition satisfies the distributive law so that the 1-coloured operad $G$ is Koszul.
In fact, this condition makes sense even if we assume that the associative product is {\em not}\,
commutative so that one might attempt to define an operad of {\em non-commutative}\, Gerstenhaber
algebras as a 1-coloured operad generated by associative non-commutative product
 product of degree 0,
$ \Ba{c}\begin{xy}
 <0mm,0.65mm>*{};<0mm,3.5mm>*{}**@{-},
 <0.5mm,-0.5mm>*{};<2.2mm,-2.2mm>*{}**@{-},
 <-0.48mm,-0.48mm>*{};<-2.2mm,-2.2mm>*{}**@{-},
 <0mm,0mm>*{\circ};<0mm,0mm>*{}**@{},
 <0.5mm,-0.5mm>*{};<2.8mm,-3.3mm>*{_2}**@{},
 <-0.48mm,-0.48mm>*{};<-2.7mm,-3.3mm>*{_1}**@{},
 \end{xy}\Ea
\neq
\Ba{c}
\begin{xy}
 <0mm,0.65mm>*{};<0mm,3.5mm>*{}**@{-},
 <0.5mm,-0.5mm>*{};<2.2mm,-2.2mm>*{}**@{-},
 <-0.48mm,-0.48mm>*{};<-2.2mm,-2.2mm>*{}**@{-},
 <0mm,0mm>*{\circ};<0mm,0mm>*{}**@{},
 <0.5mm,-0.5mm>*{};<2.8mm,-3.3mm>*{_1}**@{},
 <-0.48mm,-0.48mm>*{};<-2.7mm,-3.3mm>*{_2}**@{},
 \end{xy}\Ea$,
and Lie bracket of degree $-1$,
$\Ba{c}\begin{xy}
 <0mm,0.55mm>*{};<0mm,3.5mm>*{}**@{-},
 <0.5mm,-0.5mm>*{};<2.2mm,-2.2mm>*{}**@{-},
 <-0.48mm,-0.48mm>*{};<-2.2mm,-2.2mm>*{}**@{-},
 <0mm,0mm>*{\bu};<0mm,0mm>*{}**@{},
 <0.5mm,-0.5mm>*{};<2.8mm,-3.3mm>*{_2}**@{},
 <-0.48mm,-0.48mm>*{};<-2.7mm,-3.3mm>*{_1}**@{},
 \end{xy}\Ea
=
\Ba{c}
\begin{xy}
 <0mm,0.55mm>*{};<0mm,3.5mm>*{}**@{-},
 <0.5mm,-0.5mm>*{};<2.2mm,-2.2mm>*{}**@{-},
 <-0.48mm,-0.48mm>*{};<-2.2mm,-2.2mm>*{}**@{-},
 <0mm,0mm>*{\bu};<0mm,0mm>*{}**@{},
 <0.5mm,-0.5mm>*{};<2.8mm,-3.3mm>*{_1}**@{},
 <-0.48mm,-0.48mm>*{};<-2.7mm,-3.3mm>*{_2}**@{},
 \end{xy}\Ea$,
formally satisfying the same relations as the ones in the operad $\cG$. However, the
compatibility condition (\ref{2: Gertenshaber compatibility}) now fails to obey the
distributive law (there are new unwanted relations already for graphs with three vertices,
see Remark 1.7 in \cite{Me1}), and the resulting 1-coloured operad fails to be Koszul.
However, this problem with non-Koszulness  disappears if we think of the generating
operations as living in two different (say, dashed and straight) colours,
$$
\begin{xy}
 <0mm,0.65mm>*{};<0mm,3.5mm>*{}**@{.},
 <0.5mm,-0.5mm>*{};<2.2mm,-2.2mm>*{}**@{.},
 <-0.48mm,-0.48mm>*{};<-2.2mm,-2.2mm>*{}**@{.},
 <0mm,0mm>*{\circ};<0mm,0mm>*{}**@{},
 <0.5mm,-0.5mm>*{};<2.8mm,-3.3mm>*{_2}**@{},
 <-0.48mm,-0.48mm>*{};<-2.7mm,-3.3mm>*{_1}**@{},
 \end{xy}
\neq
\begin{xy}
 <0mm,0.65mm>*{};<0mm,3.5mm>*{}**@{.},
 <0.5mm,-0.5mm>*{};<2.2mm,-2.2mm>*{}**@{.},
 <-0.48mm,-0.48mm>*{};<-2.2mm,-2.2mm>*{}**@{.},
 <0mm,0mm>*{\circ};<0mm,0mm>*{}**@{},
 <0.5mm,-0.5mm>*{};<2.8mm,-3.3mm>*{_1}**@{},
 <-0.48mm,-0.48mm>*{};<-2.7mm,-3.3mm>*{_2}**@{},
 \end{xy}, \ \ \ \ \
 \begin{xy}
 <0mm,0.55mm>*{};<0mm,3.5mm>*{}**@{-},
 <0.5mm,-0.5mm>*{};<2.2mm,-2.2mm>*{}**@{-},
 <-0.48mm,-0.48mm>*{};<-2.2mm,-2.2mm>*{}**@{-},
 <0mm,0mm>*{\bu};<0mm,0mm>*{}**@{},
 <0.5mm,-0.5mm>*{};<2.8mm,-3.3mm>*{_2}**@{},
 <-0.48mm,-0.48mm>*{};<-2.7mm,-3.3mm>*{_1}**@{},
 \end{xy}
=
\begin{xy}
 <0mm,0.55mm>*{};<0mm,3.5mm>*{}**@{-},
 <0.5mm,-0.5mm>*{};<2.2mm,-2.2mm>*{}**@{-},
 <-0.48mm,-0.48mm>*{};<-2.2mm,-2.2mm>*{}**@{-},
 <0mm,0mm>*{\bu};<0mm,0mm>*{}**@{},
 <0.5mm,-0.5mm>*{};<2.8mm,-3.3mm>*{_1}**@{},
 <-0.48mm,-0.48mm>*{};<-2.7mm,-3.3mm>*{_2}**@{},
 \end{xy}.
$$
To make sense of the Gerstenhaber compatibility condition (\ref{2: Gertenshaber compatibility})
in two colours, we can notice that the generator
$\Ba{c}\begin{xy}
 <0mm,0.55mm>*{};<0mm,3.5mm>*{}**@{-},
 <0.5mm,-0.5mm>*{};<2.2mm,-2.2mm>*{}**@{-},
 <-0.48mm,-0.48mm>*{};<-2.2mm,-2.2mm>*{}**@{-},
 <0mm,0mm>*{\bu};<0mm,0mm>*{}**@{},
 <0.5mm,-0.5mm>*{};<2.8mm,-3.3mm>*{_2}**@{},
 <-0.48mm,-0.48mm>*{};<-2.7mm,-3.3mm>*{_1}**@{},
 \end{xy}\Ea$
 plays a two-fold role in the compatibility conditions of the operad $G$: it represents
 a Lie algebra structure, and also a {\em morphism}\ from that Lie algebra
into the Lie algebra of derivations of the associative algebra represented by
$\Ba{c}\begin{xy}
 <0mm,0.65mm>*{};<0mm,3.5mm>*{}**@{.},
 <0.5mm,-0.5mm>*{};<2.2mm,-2.2mm>*{}**@{.},
 <-0.48mm,-0.48mm>*{};<-2.2mm,-2.2mm>*{}**@{.},
 <0mm,0mm>*{\circ};<0mm,0mm>*{}**@{},
 <0.5mm,-0.5mm>*{};<2.8mm,-3.3mm>*{_2}**@{},
 <-0.48mm,-0.48mm>*{};<-2.7mm,-3.3mm>*{_1}**@{},
 \end{xy}\Ea$.
In the two coloured version we have to assign these two roles to two different actors, that is,
 we have to introduce a new degree $-1$  generator,
 $\Ba{c}\begin{xy}
 <0mm,0.55mm>*{};<0mm,3.5mm>*{}**@{.},
 <0.5mm,-0.5mm>*{};<2.2mm,-2.2mm>*{}**@{.},
 <-0.48mm,-0.48mm>*{};<-2.2mm,-2.2mm>*{}**@{-},
 <0mm,0mm>*{\bu};<0mm,0mm>*{}**@{},
 <0.5mm,-0.5mm>*{};<2.8mm,-3.3mm>*{_2}**@{},
 <-0.48mm,-0.48mm>*{};<-2.7mm,-3.3mm>*{_1}**@{},
 \end{xy}\Ea$,
for the role of the {\em morphism}, and then substitute (\ref{2: Gertenshaber compatibility})
with the following two relations,
\Beq\label{2: noncomm Gerstenhaber condition}
\Ba{c}
\begin{xy}
 <0mm,0mm>*{\bu};<0mm,0mm>*{}**@{.},
 <0mm,0.49mm>*{};<0mm,3.5mm>*{}**@{-},
 <0.49mm,-0.49mm>*{};<1.9mm,-1.9mm>*{}**@{.},
 <-0.5mm,-0.5mm>*{};<-1.9mm,-1.9mm>*{}**@{-},
 <2.3mm,-2.3mm>*{\circ};<-2.3mm,-2.3mm>*{}**@{},
 <1.8mm,-2.8mm>*{};<0mm,-4.9mm>*{}**@{.},
 <2.8mm,-2.9mm>*{};<4.6mm,-4.9mm>*{}**@{.},
   <0.49mm,-0.49mm>*{};<-2.8mm,-3.2mm>*{_1}**@{},
   <-1.8mm,-2.8mm>*{};<0.1mm,-6.2mm>*{_2}**@{},
   <-2.8mm,-2.9mm>*{};<5.3mm,-6.2mm>*{_3}**@{},
 \end{xy}
 \Ea
\ = \
\Ba{c}
\begin{xy}
 <0mm,0mm>*{\circ};<0mm,0mm>*{}**@{},
 <0mm,0.49mm>*{};<0mm,3.5mm>*{}**@{.},
 <0.49mm,-0.49mm>*{};<1.9mm,-1.9mm>*{}**@{.},
 <-0.5mm,-0.5mm>*{};<-1.9mm,-1.9mm>*{}**@{.},
 <-2.3mm,-2.3mm>*{\bu};<-2.3mm,-2.3mm>*{}**@{},
 <-1.8mm,-2.8mm>*{};<0mm,-4.9mm>*{}**@{.},
 <-2.8mm,-2.9mm>*{};<-4.6mm,-4.9mm>*{}**@{-},
   <0.49mm,-0.49mm>*{};<2.7mm,-3.2mm>*{_3}**@{},
   <-1.8mm,-2.8mm>*{};<0.4mm,-6.2mm>*{_2}**@{},
   <-2.8mm,-2.9mm>*{};<-5.1mm,-6.2mm>*{_1}**@{},
 \end{xy}\Ea
\ + \
\Ba{c}
\begin{xy}
 <0mm,0mm>*{\circ};<0mm,0mm>*{}**@{},
 <0mm,0.49mm>*{};<0mm,3.5mm>*{}**@{.},
 <0.49mm,-0.49mm>*{};<1.9mm,-1.9mm>*{}**@{.},
 <-0.5mm,-0.5mm>*{};<-1.9mm,-1.9mm>*{}**@{.},
 <2.3mm,-2.3mm>*{\bu};<-2.3mm,-2.3mm>*{}**@{},
 <1.8mm,-2.8mm>*{};<0mm,-4.9mm>*{}**@{-},
 <2.8mm,-2.9mm>*{};<4.6mm,-4.9mm>*{}**@{.},
   <0.49mm,-0.49mm>*{};<-2.8mm,-3.2mm>*{_2}**@{},
   <-1.8mm,-2.8mm>*{};<0.1mm,-6.2mm>*{_1}**@{},
   <-2.8mm,-2.9mm>*{};<5.3mm,-6.2mm>*{_3}**@{},
 \end{xy}
 \Ea,
 \ \ \ \ \ \ \
 \Ba{c}
\begin{xy}
 <0mm,0mm>*{\bu};<0mm,0mm>*{}**@{},
 <0mm,0.49mm>*{};<0mm,3.5mm>*{}**@{.},
 <0.49mm,-0.49mm>*{};<1.9mm,-1.9mm>*{}**@{.},
 <-0.5mm,-0.5mm>*{};<-1.9mm,-1.9mm>*{}**@{-},
 <-2.3mm,-2.3mm>*{\bu};<-2.3mm,-2.3mm>*{}**@{},
 <-1.8mm,-2.8mm>*{};<0mm,-4.9mm>*{}**@{-},
 <-2.8mm,-2.9mm>*{};<-4.6mm,-4.9mm>*{}**@{-},
   <0.49mm,-0.49mm>*{};<2.7mm,-3.2mm>*{_3}**@{},
   <-1.8mm,-2.8mm>*{};<0.4mm,-6.2mm>*{_2}**@{},
   <-2.8mm,-2.9mm>*{};<-5.1mm,-6.2mm>*{_1}**@{},
 \end{xy}\Ea
\ = \
\Ba{c}
\begin{xy}
 <0mm,0mm>*{\bu};<0mm,0mm>*{}**@{},
 <0mm,0.49mm>*{};<0mm,3.5mm>*{}**@{.},
 <0.49mm,-0.49mm>*{};<1.9mm,-1.9mm>*{}**@{.},
 <-0.5mm,-0.5mm>*{};<-1.9mm,-1.9mm>*{}**@{-},
 <2.3mm,-2.3mm>*{\bu};<-2.3mm,-2.3mm>*{}**@{},
 <1.8mm,-2.8mm>*{};<0mm,-4.9mm>*{}**@{-},
 <2.8mm,-2.9mm>*{};<4.6mm,-4.9mm>*{}**@{.},
   <0.49mm,-0.49mm>*{};<-2.8mm,-3.2mm>*{_1}**@{},
   <-1.8mm,-2.8mm>*{};<0.1mm,-6.2mm>*{_2}**@{},
   <-2.8mm,-2.9mm>*{};<5.3mm,-6.2mm>*{_3}**@{},
 \end{xy}
 \Ea
 \ + \
\Ba{c}
\begin{xy}
 <0mm,0mm>*{\bu};<0mm,0mm>*{}**@{},
 <0mm,0.49mm>*{};<0mm,3.5mm>*{}**@{.},
 <0.49mm,-0.49mm>*{};<1.9mm,-1.9mm>*{}**@{.},
 <-0.5mm,-0.5mm>*{};<-1.9mm,-1.9mm>*{}**@{-},
 <2.3mm,-2.3mm>*{\bu};<-2.3mm,-2.3mm>*{}**@{},
 <1.8mm,-2.8mm>*{};<0mm,-4.9mm>*{}**@{-},
 <2.8mm,-2.9mm>*{};<4.6mm,-4.9mm>*{}**@{.},
   <0.49mm,-0.49mm>*{};<-2.8mm,-3.2mm>*{_2}**@{},
   <-1.8mm,-2.8mm>*{};<0.1mm,-6.2mm>*{_1}**@{},
   <-2.8mm,-2.9mm>*{};<5.3mm,-6.2mm>*{_3}**@{},
 \end{xy}
 \Ea
\Eeq

The 2-coloured operad generated by binary operations $
\begin{xy}
 <0mm,0.65mm>*{};<0mm,3.5mm>*{}**@{.},
 <0.5mm,-0.5mm>*{};<2.2mm,-2.2mm>*{}**@{.},
 <-0.48mm,-0.48mm>*{};<-2.2mm,-2.2mm>*{}**@{.},
 <0mm,0mm>*{\circ};<0mm,0mm>*{}**@{},
 <0.5mm,-0.5mm>*{};<2.8mm,-3.3mm>*{_2}**@{},
 <-0.48mm,-0.48mm>*{};<-2.7mm,-3.3mm>*{_1}**@{},
 \end{xy}
\neq
\begin{xy}
 <0mm,0.65mm>*{};<0mm,3.5mm>*{}**@{.},
 <0.5mm,-0.5mm>*{};<2.2mm,-2.2mm>*{}**@{.},
 <-0.48mm,-0.48mm>*{};<-2.2mm,-2.2mm>*{}**@{.},
 <0mm,0mm>*{\circ};<0mm,0mm>*{}**@{},
 <0.5mm,-0.5mm>*{};<2.8mm,-3.3mm>*{_1}**@{},
 <-0.48mm,-0.48mm>*{};<-2.7mm,-3.3mm>*{_2}**@{},
 \end{xy}$,
 $\begin{xy}
 <0mm,0.55mm>*{};<0mm,3.5mm>*{}**@{-},
 <0.5mm,-0.5mm>*{};<2.2mm,-2.2mm>*{}**@{-},
 <-0.48mm,-0.48mm>*{};<-2.2mm,-2.2mm>*{}**@{-},
 <0mm,0mm>*{\bu};<0mm,0mm>*{}**@{},
 <0.5mm,-0.5mm>*{};<2.8mm,-3.3mm>*{_2}**@{},
 <-0.48mm,-0.48mm>*{};<-2.7mm,-3.3mm>*{_1}**@{},
 \end{xy}
=
\begin{xy}
 <0mm,0.55mm>*{};<0mm,3.5mm>*{}**@{-},
 <0.5mm,-0.5mm>*{};<2.2mm,-2.2mm>*{}**@{-},
 <-0.48mm,-0.48mm>*{};<-2.2mm,-2.2mm>*{}**@{-},
 <0mm,0mm>*{\bu};<0mm,0mm>*{}**@{},
 <0.5mm,-0.5mm>*{};<2.8mm,-3.3mm>*{_1}**@{},
 <-0.48mm,-0.48mm>*{};<-2.7mm,-3.3mm>*{_2}**@{},
 \end{xy}$
 and
 $\Ba{c}\begin{xy}
 <0mm,0.55mm>*{};<0mm,3.5mm>*{}**@{.},
 <0.5mm,-0.5mm>*{};<2.2mm,-2.2mm>*{}**@{.},
 <-0.48mm,-0.48mm>*{};<-2.2mm,-2.2mm>*{}**@{-},
 <0mm,0mm>*{\bu};<0mm,0mm>*{}**@{},
 <0.5mm,-0.5mm>*{};<2.8mm,-3.3mm>*{_2}**@{},
 <-0.48mm,-0.48mm>*{};<-2.7mm,-3.3mm>*{_1}**@{},
 \end{xy}\Ea$,
\Beq\label{2: assoc relations}
\Ba{c}
\begin{xy}
 <0mm,0mm>*{\circ};<0mm,0mm>*{}**@{},
 <0mm,-0.49mm>*{};<0mm,3.5mm>*{}**@{.},
 <0.49mm,-0.49mm>*{};<1.9mm,-1.9mm>*{}**@{.},
 <-0.5mm,-0.5mm>*{};<-1.9mm,-1.9mm>*{}**@{.},
 <2.3mm,-2.3mm>*{\circ};<2.3mm,-2.3mm>*{}**@{},
 <1.8mm,-2.8mm>*{};<0mm,-4.9mm>*{}**@{.},
 <2.8mm,-2.9mm>*{};<4.6mm,-4.9mm>*{}**@{.},
   <0.49mm,-0.49mm>*{};<-2.7mm,-3.2mm>*{_1}**@{},
   <-1.8mm,-2.8mm>*{};<0.4mm,-6.2mm>*{_2}**@{},
   <-2.8mm,-2.9mm>*{};<5.1mm,-6.2mm>*{_3}**@{},
 \end{xy}\Ea
\ = \
\Ba{c}
\begin{xy}
 <0mm,0mm>*{\circ};<0mm,0mm>*{}**@{},
 <0mm,0.49mm>*{};<0mm,3.5mm>*{}**@{.},
 <0.49mm,-0.49mm>*{};<1.9mm,-1.9mm>*{}**@{.},
 <-0.5mm,-0.5mm>*{};<-1.9mm,-1.9mm>*{}**@{.},
 <-2.3mm,-2.3mm>*{\circ};<-2.3mm,-2.3mm>*{}**@{},
 <-1.8mm,-2.8mm>*{};<0mm,-4.9mm>*{}**@{.},
 <-2.8mm,-2.9mm>*{};<-4.6mm,-4.9mm>*{}**@{.},
   <0.49mm,-0.49mm>*{};<2.7mm,-3.2mm>*{_3}**@{},
   <-1.8mm,-2.8mm>*{};<0.4mm,-6.2mm>*{_2}**@{},
   <-2.8mm,-2.9mm>*{};<-5.1mm,-6.2mm>*{_1}**@{},
 \end{xy}\Ea,
 \Eeq
  Jacobi relations for the Lie brackets,
\Beq\label{2: Jacobi relations}
\Ba{c}
\begin{xy}
 <0mm,0mm>*{\bu};<0mm,0mm>*{}**@{},
 <0mm,0.49mm>*{};<0mm,3.5mm>*{}**@{-},
 <0.49mm,-0.49mm>*{};<1.9mm,-1.9mm>*{}**@{-},
 <-0.5mm,-0.5mm>*{};<-1.9mm,-1.9mm>*{}**@{-},
 <-2.3mm,-2.3mm>*{\bu};<-2.3mm,-2.3mm>*{}**@{},
 <-1.8mm,-2.8mm>*{};<0mm,-4.9mm>*{}**@{-},
 <-2.8mm,-2.9mm>*{};<-4.6mm,-4.9mm>*{}**@{-},
   <0.49mm,-0.49mm>*{};<2.7mm,-3.2mm>*{_3}**@{},
   <-1.8mm,-2.8mm>*{};<0.4mm,-6.2mm>*{_2}**@{},
   <-2.8mm,-2.9mm>*{};<-5.1mm,-6.2mm>*{_1}**@{},
 \end{xy}\Ea
\ + \
\Ba{c}\begin{xy}
 <0mm,0mm>*{\bu};<0mm,0mm>*{}**@{},
 <0mm,0.49mm>*{};<0mm,3.5mm>*{}**@{-},
 <0.49mm,-0.49mm>*{};<1.9mm,-1.9mm>*{}**@{-},
 <-0.5mm,-0.5mm>*{};<-1.9mm,-1.9mm>*{}**@{-},
 <-2.3mm,-2.3mm>*{\bu};<-2.3mm,-2.3mm>*{}**@{},
 <-1.8mm,-2.8mm>*{};<0mm,-4.9mm>*{}**@{-},
 <-2.8mm,-2.9mm>*{};<-4.6mm,-4.9mm>*{}**@{-},
   <0.49mm,-0.49mm>*{};<2.7mm,-3.2mm>*{_2}**@{},
   <-1.8mm,-2.8mm>*{};<0.4mm,-6.2mm>*{_1}**@{},
   <-2.8mm,-2.9mm>*{};<-5.1mm,-6.2mm>*{_3}**@{},
 \end{xy}\Ea
\ + \
\Ba{c}
\begin{xy}
 <0mm,0mm>*{\bu};<0mm,0mm>*{}**@{},
 <0mm,0.49mm>*{};<0mm,3.5mm>*{}**@{-},
 <0.49mm,-0.49mm>*{};<1.9mm,-1.9mm>*{}**@{-},
 <-0.5mm,-0.5mm>*{};<-1.9mm,-1.9mm>*{}**@{-},
 <-2.3mm,-2.3mm>*{\bu};<-2.3mm,-2.3mm>*{}**@{},
 <-1.8mm,-2.8mm>*{};<0mm,-4.9mm>*{}**@{-},
 <-2.8mm,-2.9mm>*{};<-4.6mm,-4.9mm>*{}**@{-},
   <0.49mm,-0.49mm>*{};<2.7mm,-3.2mm>*{_1}**@{},
   <-1.8mm,-2.8mm>*{};<0.4mm,-6.2mm>*{_3}**@{},
   <-2.8mm,-2.9mm>*{};<-5.1mm,-6.2mm>*{_2}**@{},
 \end{xy}\Ea
\ = \ 0,
\Eeq
 and the compatibility  relations (\ref{2: noncomm Gerstenhaber condition})
 was introduced in \cite{FGV} (with slightly different grading conventions which in two colours
 are irrelevant) under the name of the {\em operad of Leibniz pairs}. However algebras over the operad
 of Leibniz pairs have nothing to do with pairs of Leibniz algebras. We prefer to call this
 quadratic operad the {\em 2-coloured operad of noncommutative Gerstenhaber algebras}\,
 ($nc\cG$ for short) as this name specifies its structure non-ambiguously; this is the only natural
  way to generalize the notion of Gerstenhaber algebras to the case of a non-commutative product
  while keeping the Koszulness property. Moreover, any Gerstenhaber algebra is automatically an
  algebra over $nc\cG$. In particular, for any smooth manifold $M$ the associated space of
  polyvector fields, $\cT_{poly}(M)$ equipped with the Schouten bracket $[\ ,\ ]_S$ and
  the wedge product $\wedge$ is an $nc\cG$-algebra. It was proven
 in \cite{Ta1} that
 $(\cT_{poly}(\R^d), [\ ,\ ]_S, \wedge)$
 is rigid as a  $\cG_\infty$ algebra. It follows from Willwacher's proof \cite{Wi} of the Furusho theorem that
 $(\cT_{poly}(\R^d), [\ ,\ ]_S, \wedge)$ admits a unique (up to homotopy and rescalings)
 universal $nc\cG_\infty$ deformation whose explicit structure is described in \cite{A} (see also (\ref{4: Exotic ncG}) below  for its explicit graph representation).


\subsection{Configuration space model for the 4-coloured operad of morphisms
of $nc\cG_\infty$-algebras}\label{2: subsection on Mor(NCG)} A geometric model
for the 4-coloured operad
of morphisms of OCHA algebras was given in \cite{Me3}. The same ideas work for the operad,
$\cM or(nc\cG_\infty)$, of morphisms
of $nc\cG_\infty$-algebras provided one replaces everywhere in \S 6 of   \cite{Me3} the upper-plane  $\bbH$ with the full complex plane $\C$.

\bip

\bip
{\large
\section{\bf T.\ Willwacher's theorems}
}

\bip
\subsection{Universal deformations of the Schouten bracket} The deformation complex
of the graded Lie algebra $(\cT_{poly}(\R^d), [\ ,\ ]_S)$ is the graded Lie algebra,
$$
\CoDer\left(\underbrace{\odot^\bu(\cT_{poly}(\R^d)[2])}_{\mathrm{standard\ coalgebra}\atop \mathrm{structure}}\right) =\prod_{n\geq 0} \Hom(\odot^n \cT_{poly}(\R^d),
\cT_{poly}(\R^d))[2-2n]
$$
 of coderivations of the graded-cocommutative
coalgebra $\odot^\bu(\cT_{poly}(\R^d)[2])$ equipped with the differential, $\delta$, given by
$$
\delta (D):=  [\ ,\ ]_S \circ D - (-1)^{|D|} D\circ  [\ ,\ ]_S,\ \ \forall\
D\in \CoDer\left(\odot^\bu(\cT_{poly}(\R^d)[2])\right).
$$
Here $\circ$ stands for the composition of coderivations.
There is a universal (i.e.\ independent of the dimension $d$) version of this deformation complex,
 $\sG\sC_2$, which was introduced by Kontsevich in \cite{Ko1} and studied in detail in
 \cite{Wi}. In this subsection we recall some ideas, results and notations of \cite{Wi} which we later use to prove our main theorem.

\subsection{Operad $\cG ra$}\label{subsection 2: Gra}
To define Kontsevich's dg Lie algebra $\sG\sC_2$ it is easiest to start by defining a certain
operad of graphs. For arbitrary integers $n\geq 1$ and $l\geq 0$ let ${\sG}_{n,l}$ stand for the
set of graphs\ $\{\Ga\}$ with $n$ vertices and $l$ edges
such that (i) the vertices of $\Ga$ are labelled by elements of $[n]:=\{1,\ldots, n\}$,
(ii) the set of edges, $E(\Ga)$, is totally ordered up to an even permutation (that is, {\em oriented}); it has at most two different orientations. For $\Ga\in {\sG}_{n,l}$ we denote by
$\Ga_{opp}$ the oppositely oriented graph.
Let $\K\langle \sG_{n,l}\rangle$  be the vector space over  a field $\K$ spanned by isomorphism
classes, $[\Ga]$, of elements of $\sG_{n,l}$ modulo the
relation\footnote{Abusing notations we identify from now an equivalence class $[\Ga]$ with any
of its representative $\Ga$.}  $\Ga_{opp}=-\Ga$, and consider the $\Z$-graded $\bS_n$-module,
$$
\cG ra (n):=\bigoplus_{l=0}^\infty \K\langle \sG_{n,l}\rangle[l].
$$
For example,  $\xy
(0,2)*{_{1}},
(7,2)*{_{2}},
 (0,0)*{\bullet}="a",
(7,0)*{\bu}="b",
\ar @{-} "a";"b" <0pt>
\endxy$ is a degree $-1$ element in  $\cG ra(2)$.
 The $\bS$-module, $\cG ra :=\{\cG ra (n)\}_{n\geq 1}$, is naturally an operad with the
 operadic compositions given by
\Beq\label{3: operad comp in Gra}
\Ba{rccc}
\circ_i: & \cG ra (n)\ot \cG ra (m) & \lon &  \cG ra (m+n-1)\\
&  \Ga_1 \ot \Ga_2   &\lon & \sum_{\Ga\in \sG_{\Ga_1, \Ga_2}^i} (-1)^{\sigma_\Ga} \Ga
\Ea
\Eeq
where $ \sG_{\Ga_1, \Ga_2}^i$ is the subset of $\sG_{n+m-1, \# E(\Ga_1) + \#E(\Ga_2)}$ consisting
of graphs, $\Ga$, satisfying the condition: the full subgraph of $\Ga$ spanned by the vertices labeled by
the set $\{i,i+1, \ldots, i+m-1\}$ is isomorphic to $\Ga_2$ and the quotient graph, $\Ga/\Ga_2$,
obtained by contracting that subgraph to a single vertex, is isomorphic to $\Ga_1$ (see \S 2 in
\cite{Wi} or \S 7 in \cite{Me3} for examples).
The sign $(-1)^{\sigma_\Gamma}$ is determined by the equality
$
\wedge_{e\in E(\Ga)}e= (-1)^{\sigma_\Gamma}(\wedge_{e'\in E(\Ga_1)}e')
\wedge (\bigwedge_{e''\in E(\Ga_2)}e'')$ where the edge products over the sets
$E(\Ga_1)$ and $E(\Ga_1)$ are taken in accordance with the given orientations.
The unique element in $\sG_{1,0}$ serves as the unit element in the operad  $\cG ra$.


\subsection{A canonical representation of $\cG ra$ in $\cT_{poly}(\R^d)$}\label{3: subsec repr of Gra}
The operad  $\cG ra$ has a natural representation in the vector space
$\cT_{poly}(\R^d)[2]$ for any dimension $d$,
\Beq\label{3: Gra representation}
\Ba{rccc}
\rho: & \cG ra(n) & \lon & \cE  nd_{\cT_{poly}(\R^d)}(n)=
\Hom( \cT_{poly}(\R^d)^{\ot n},\cT_{poly}(\R^d))\\
      & \Ga &\lon & \Phi_\Ga
\Ea
\Eeq
given by the formula,
$$
\Phi_\Ga(\ga_1,\ldots, \ga_n) :=\mu\left(\prod_{e\in E(\Ga)}
\Delta_e \left(\ga_1(x_{(1)}, \psi_{(1)})\ot \ga_2(x_{(2)}, \psi_{(2)})\ot \ldots\ot
\ga_n(x_{(n)}, \psi_{(n)}) \right)\right)
$$
where, for an edge $e$ connecting vertices labeled by integers $i$ and $j$,
$$
\Delta_e= \sum_{a=1}^n \frac{\p}{\p x_{(i)}^a} \ot  \frac{\p}{\p \psi_{(j)a}}
+
\frac{\p}{\psi_{(i)a}}\ot  \frac{\p}{\p x_{(j)}^a}
$$
and $\mu$ is the multiplication map,
$$
\Ba{rccc}
\mu:&   \cT_{poly}(\R^d)^{\ot n} & \lon & \cT_{poly}(\R^d)\\
   & \ga_1\ot \ga_2\ot \ldots \ot \ga_n &\lon & \ga_1\wedge \ga_2\wedge \cdots
   \wedge \ga_n.
\Ea
$$
Here we used a coordinate identification,
$\cT_{poly}(\R^d)=C^\infty(x^1,\ldots, x^d)[\psi_1, \ldots, \psi_d]$, where
$C^\infty(x^1,\ldots, x^d)$ is the ring of smooth functions of coordinates $x^1,\ldots, x^d$
on $\R^d$, and $\psi_a$ are formal variables of degree one symbolizing $\p/\p x^a$.

\subsection{Kontsevich graph complex} There is a morphism of operads \cite{Wi3}
$$
\cG \lon \cG ra
$$
given on the generators of the operad of Gerstenhaber algebras by
\Beq\label{2: Comm to Gra}
\Ba{c}
\xy
 <0mm,0.55mm>*{};<0mm,3.5mm>*{}**@{-},
 <0.5mm,-0.5mm>*{};<2.2mm,-2.2mm>*{}**@{-},
 <-0.48mm,-0.48mm>*{};<-2.2mm,-2.2mm>*{}**@{-},
 <0mm,0mm>*{\circ};<0mm,0mm>*{}**@{},
 <0.5mm,-0.5mm>*{};<2.7mm,-3.2mm>*{_2}**@{},
 <-0.48mm,-0.48mm>*{};<-2.7mm,-3.2mm>*{_1}**@{},
 \endxy\Ea \ \ \lon \ \ \ \xy
(0,2)*{_{1}},
(5,2)*{_{2}},
 (0,0)*{\bullet}="a",
(5,0)*{\bu}="b",
\endxy
\Eeq
\Beq\label{2: Lie to Gra}
\Ba{c}
\xy
 <0mm,0.55mm>*{};<0mm,3.5mm>*{}**@{-},
 <0.5mm,-0.5mm>*{};<2.2mm,-2.2mm>*{}**@{-},
 <-0.48mm,-0.48mm>*{};<-2.2mm,-2.2mm>*{}**@{-},
 <0mm,0mm>*{\bu};<0mm,0mm>*{}**@{},
 <0.5mm,-0.5mm>*{};<2.7mm,-3.2mm>*{_2}**@{},
 <-0.48mm,-0.48mm>*{};<-2.7mm,-3.2mm>*{_1}**@{},
 \endxy\Ea
   \ \ \ \ \lon \ \ \ \ \xy
(0,2)*{_{1}},
(5,2)*{_{2}},
 (0,0)*{\bullet}="a",
(5,0)*{\bu}="b",
\ar @{-} "a";"b" <0pt>
\endxy
\Eeq
The latter map also gives us a canonical morphism of operads
$$
i: \caL ie\{1\} \lon \cG ra.
$$
The {\em full Kontsevich graph complex}\, $\mathsf{fGC}_2$ is, by definition, the deformation
complex controlling deformations of the morphism $i$,
$$
\mathsf{fGC}_2:= \Def(\caL ie\{1\} \rar \cG ra)
$$
 There are several explicit constructions of deformation complexes of (pr)operadic morphisms
  given, for example, in \cite{MV}. To construct  $\Def(\caL ie\{1\} \rar \cG ra)$ one has to
  replace $\caL ie\{1\}$ by its the minimal resolution, $\caL ie\{1\}_\infty$,  which is
  a quasi-free dg operad generated by the $\bS$-module
 $$
 E=\left\{E(n):=\id_n[2n-3]\right\}.
 $$
Then, as a $\Z$-graded vector space,
$$
\Def(\caL ie\{1\} \rar \cG ra)\equiv \Def(\caL ie\{1\}_\infty \rar \cG ra):= \bigoplus_{n\geq 0}\Hom_{\bS_n}(E(n), Gra(n))[-1]=
 \bigoplus_{n\geq 0} Gra(n)^{\bS_n}[2-2n],
$$
i.e.\ an element of $\mathsf{fGC}_2$
can be understood as an $\bS_n$-symmetrization a  of graph from  $\sG_{n,l}$  to which we assign  the degree $2n-l-2$, for example
$$
\Ba{c}\xy
(0,2)*{_{1}},
(6,2)*{_{2}},
 (0,0)*{\bullet}="a",
(6,0)*{\bu}="b",
\ar @{-} "a";"b" <0pt>
\endxy \Ea + \Ba{c}\xy
(0,2)*{_{2}},
(6,2)*{_{1}},
 (0,0)*{\bullet}="a",
(6,0)*{\bu}="b",
\ar @{-} "a";"b" <0pt>
\endxy\Ea=: \xy
 (0,0)*{\bullet}="a",
(6,0)*{\bu}="b",
\ar @{-} "a";"b" <0pt>
\endxy
$$
is a degree 1 element in $\mathsf{fGC}_2$. As labelling of vertices of elements from $\mathsf{fGC}_2$ by integers is symmetrized, we often represent such elements as a single graph with vertices {\em unlabelled}, e.g.
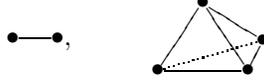
$$
 \xy
 (0,0)*{\bullet}="a",
(6,0)*{\bu}="b",
\ar @{-} "a";"b" <0pt>
\endxy, \ \ \ \ \ \ \ \Ba{c}
\xy
 (0,0)*{\bullet}="a",
(12,0)*{\bullet}="b",
(6,9)*{\bullet}="c",
(14,4)*{\bullet}="d",
\ar @{-} "a";"b" <0pt>
\ar @{-} "a";"c" <0pt>
\ar @{-} "b";"c" <0pt>
\ar @{-} "d";"c" <0pt>
\ar @{-} "b";"d" <0pt>
\ar @{.} "d";"a" <0pt>
\endxy
\Ea
$$
One should not forget, however, that such a graph is in reality a symmetrization
sum of some labelled graph from  $\sG_{n,l}$.

\sip

The Lie algebra structure in $\mathsf{fGC}_2= \Def(\caL ie\{1\}_\infty \rar \cG ra)$ is completely determined by the differential on $\caL ie\{1\}_\infty$ \cite{MV}. It is an elementary exercise
to see that the Lie brackets in $\mathsf{fGC}_2$ can expressed in terms of operadic composition
in $\cG ra$ as follows,
$$
[\Ga, \Ga']:= Sym( \Ga\circ_1 \Ga' - (-1)^{|\Ga||\Ga'|}\Ga\circ_1 \Ga),
$$
where $Sym$ stands for the symmetrization of vertex labels. The usefulness of this Lie algebra structure on  $\mathsf{fGC}_2:= \Def(\caL ie\{1\}_\infty \rar \cG ra)$ stems from the fact \cite{MV} that the set of its Maurer-Cartan elements is in one-to-one correspondence
with morphisms of operads  $\caL ie\{1\}_\infty \rar \cG ra$. It is easy to check that
the element  $\xy
 (0,0)*{\bullet}="a",
(5,0)*{\bu}="b",
\ar @{-} "a";"b" <0pt>
\endxy$ is  Maurer-Cartan,
$$
[ \xy
 (0,0)*{\bullet}="a",
(5,0)*{\bu}="b",
\ar @{-} "a";"b" <0pt>
\endxy, \xy
 (0,0)*{\bullet}="a",
(5,0)*{\bu}="b",
\ar @{-} "a";"b" <0pt>
\endxy]=0.
$$
It corresponds precisely to the morphism (\ref{2: Lie to Gra}). This element makes
 $\mathsf{fGC}_2$ into a complex with the differential
$$
\delta_\bubu:=[\xy
 (0,0)*{\bullet}="a",
(5,0)*{\bu}="b",
\ar @{-} "a";"b" <0pt>
\endxy,\    ].
$$
 This dg Lie
algebra contains a dg Lie subalgebra, $\sG\sC_2$, spanned by connected graphs with at
least trivalent vertices and no tadpoles;
this subalgebra is precisely
the original (odd) {\em Kontsevich graph complex} \cite{Ko1, Wi}.
One of the main theorems of \cite{Wi} asserts an isomorphism of Lie algebras,
$$
H^0(\sG\sC_2, \delta_\bubu)\simeq \fg\fr\ft,
$$
where $\fg\fr\ft$ stands for the Grothendieck-Teichm\"uller Lie algebra and $H^0$ for
cohomology in degree zero.

\sip

Note that the canonical representation (\ref{3: Gra representation}) induces a morphism
of dg Lie algebras,
$$
\rho^{ind}: \ssf\sG\sC_2= \Def\left(\caL ie\{1\}_\infty \rar \cG ra\right)
 \lon   \Def\left(\caL ie\{1\}_\infty \rar   \cE  nd_{\cT_{poly}(\R^d)}\right)  =   \CoDer\left(\odot^\bu(\cT_{poly}(\R^d)[2])\right).
$$
The image of this map consists of coderivations of the coalgebra $\odot^\bu(\cT_{poly}(\R^d)[2]$ which are {\em universal}\, i.e.\ make sense in any dimension. In particular, $\rho(\xy
 (0,0)*{\bullet}="a",
(5,0)*{\bu}="b",
\ar @{-} "a";"b" <0pt>
\endxy)$ is precisely the Schouten bracket in  $\cT_{poly}(\R^d)$. Therefore, one can say that
 the graph complex  $\ssf\sG\sC_2$ (or $\sG\sC_2$) describes {\em universal}\,
deformations of the Schouten bracket. T.\ Willwacher's theorem gives us universal homotopy actions of the Grothendieck-Teichm\"uller group $GRT=\exp(\fg\fr\ft)$ on  $\cT_{poly}(\R^d)$
by $\caL ie_\infty$ automorphisms of the Schouten bracket.

\subsection{T.\ Willwacher's twisted operad $f\cG raphs^\circlearrowleft$} For any operad $\cP$ and morphism of operads,
$\caL ie\{k\}_\infty \rar \cP$, there is an associated operad $Tw(\cP)$ whose representations,
$\rho^{tw}: Tw(\cP)\rar \cE nd_V$, can be obtained from representations, $\rho: \cP\rar \cE nd_V$, of $\cP$ by ``twisting" $\rho$
 by Maurer-Cartan elements of the associated (via the map $\caL ie\{k\}_\infty \rar \cP$) $\caL ie\{k\}_\infty$ structure on $V$. Omitting general construction (see \cite{Wi} for its details), we shall describe explicitly the dg operad
 $f\cG raphs^\circlearrowleft:=Tw(\cG ra)$ obtained from $\cG ra$ by twisting the morphism (\ref{2: Lie to Gra}).
For arbitrary integers $m\geq 1$, $n\geq 0$ and $l\geq 0$ we denote by  ${\sG}_{m,n;l}$  a set
of graphs\ $\{\Ga\}$ with $m$ white vertices, $n$ black
vertices
are and $l$ edges such that (i) the white vertices of $\Ga$ are labelled by elements of $[m]$, (ii) the black vertices of $\Ga$ are at least trivalent and are labelled by elements of $[\bar{n}]=\{\bar{1},\dots, \bar{n}\}$,
(iii) and the set of edges, $E(\Ga)$, is totally ordered up to an even permutation.
The set of black
(respectively, white) vertices of $\Ga$ will be denoted by $V_\bu(\Ga)$ (resp.\, $V_\circ(\Ga)$).

Let $\K\langle \sG_{m,n;l}\rangle$  be the vector space over  a field $\K$ spanned by isomorphism classes,
 $[\Ga]$, of elements of $\sG_{m,n;l}$ modulo the relation $\Ga_{opp}=-\Ga$, and consider
the $\Z$-graded $\bS_m$-module,
$$
f\cG raphs^\circlearrowleft (m):=\bigoplus_{l=0}^\infty\bigoplus_{n=0}^\infty \K\langle
\sG_{m,n;l}\rangle^{\bS_n}[l],
$$
where invariants are taken with respect to the permutations of $[\bar{n}]$-labellings
of black vertices. For example,  $\xy (0,2)*{_{1}},
 (0,0)*{\circ}="a",
(5,0)*{\bu}="b",
\ar @{-} "a";"b" <0pt>\endxy$ is a degree $-1$ element in  $f\cG raphs^\circlearrowleft(1)$ and $\xy
(0,2)*{_{1}},
(5,2)*{_{2}},
 (0,0)*{\circ}="a",
(5,0)*{\circ}="b", \endxy$ is a degree $0$ element in  $f\cG raphs^\circlearrowleft(2)$.
The operadic composition, $\Ga\circ_i \Ga'$, in
$$
f\cG raphs^\circlearrowleft=\left\{f\cG raphs^\circlearrowleft (m)\right\}
$$
is defined by substitution of the graph $\Ga'\in\K\langle \sG_{m',n';l}\rangle^{\bS_{n'}}$ into the $i$-th white
vertex $v$ of $\Ga\in\K\langle \sG_{m,n;l}\rangle^{\bS_{n}}$, reconnecting all edges of \(\Ga\) incident to $v$
in all possible ways to vertices of $\Ga'$ (in a full analogy to the case of $\cG ra$),
and finally symmetrizing over labellings of the $n+n'$ black vertices. Consider linear maps,
$$
\delta_\bubu \Ga := -(-1)^{|\Ga|} Sym  \left(\Ga \circ_{\bar{1}} \xy
 (0,0)*{\bullet}="a",
(5,0)*{\bu}="b",
\ar @{-} "a";"b" <0pt>
\endxy\right)
$$
and
$$
\delta_\wibu \Ga :=  Sym\left(
 \xy(0,2)*{_{1}},
 (0,0)*{\circ}="a",
(5,0)*{\bu}="b",
\ar @{-} "a";"b" <0pt>
\endxy \circ_1 \Ga\ \
- \ \ (-1)^{|\Ga|} \sum_{v\in V(\circ)} \Ga \circ_v
\xy(0,2)*{_{1}},
 (0,0)*{\circ}="a",
(5,0)*{\bu}="b",
\ar @{-} "a";"b" <0pt>
\endxy   \right)
$$
where $Sym$ stands for the symmetrization of black vertex labellings. Note that
in this case $\delta_\bubu^2\neq 0$ and $\delta_\wibu^2\neq 0$ in general, but their sum
$\delta_\wibu+\delta_\bubu$
 makes $f\cG raphs^\circlearrowleft$ into an operad of {\em complexes} \cite{Wi}.

\sip

The dg suboperad of  $f\cG raphs^\circlearrowleft$ consisting of graphs $\Ga$ which have no connected component consisting solely of black vertices is denoted in \cite{Wi} by $\cG raphs^\circlearrowleft$.  Without loss
 of much generality \cite{Wi} we may replace operads  $f\cG raphs^\circlearrowleft$ and  $\cG raphs^\circlearrowleft$ by their suboperads, $f\cG raphs$ and  respectively $\cG raphs$, consisting of graphs with no tadpoles.

\sip

There is a morphism of dg operads \cite{Wi2}
$$
\cA ss_\infty \lon \cA ss \lon \cG raphs
$$
where the first arrow is a natural projection and the second map is
given on the generators of the operad $\cA ss$ by
$$
\xy
 <0mm,0.55mm>*{};<0mm,3.5mm>*{}**@{-},
 <0.5mm,-0.5mm>*{};<2.2mm,-2.2mm>*{}**@{-},
 <-0.48mm,-0.48mm>*{};<-2.2mm,-2.2mm>*{}**@{-},
 <0mm,0mm>*{\circ};<0mm,0mm>*{}**@{},
 <0.5mm,-0.5mm>*{};<2.7mm,-3.2mm>*{_2}**@{},
 <-0.48mm,-0.48mm>*{};<-2.7mm,-3.2mm>*{_1}**@{},
 \endxy \ \ \lon \ \ \ \xy
(0,2)*{_{1}},
(5,2)*{_{2}},
 (0,0)*{\circ}="a",
(5,0)*{\circ}="b",
\endxy
$$
The standard construction \cite{MV} gives us a dg Lie algebra,
$\Def(\cA ss \rar \cG raphs)$,
whose elements, $\Ga$, are linear combinations of graphs from $\K\langle
\sG_{m,n;l}\rangle^{\bS_n}$, $m,n,l\geq 0$,
 equipped with a total order on the set of white vertices of $\Ga$
(so that in pictures we can depict vertices of such graphs as lying on a line) and with degree $2n+m-l-1$. The differential on $\Def(\cA ss_\infty \rar f\cG raphs)$ is a sum,
\Beq\label{3: differential in Def(Ass--Graphs)}
\delta= \delta_\wiwi + \delta_\wibu + \delta_\bubu,
\Eeq
where
$$
\delta_\wiwi \Ga:=
 \left( \xy
(0,2)*{_{1}},(5,2)*{_{2}},
 (0,0)*{\circ}="a",
(5,0)*{\circ}="b",
\ar @{-} "a";"b" <0pt>
\endxy \circ_1 \Ga
+
\xy
(0,2)*{_{1}},(5,2)*{_{2}},
 (0,0)*{\circ}="a",
(5,0)*{\circ}="b",
\ar @{-} "a";"b" <0pt>
\endxy \circ_2 \Ga
\right)
\  - \  (-1)^{|\Ga|}\sum_{v\in V(\circ)} \Ga \circ_v
\xy
(0,2)*{_{1}},
(5,2)*{_{2}},
 (0,0)*{\circ}="a",
(5,0)*{\circ}="b",
\ar @{-} "a";"b" <0pt>
\endxy
$$
The first cohomology group of this deformation complex was computed  in \cite{Wi,Wi2},
\Beq\label{2: H of Def(Ass-Graphs)}
H^i\left(\Def(\cA ss_\infty\hspace{-0.8mm} \rar \hspace{-0.8mm}\cG raphs)\right)=\left\{\Ba{ll}
\fg\fr\ft \ \oplus \ \R[-1] & \mathrm{for}\ i=1,\\
\R[\bS_2]                          & \mathrm{for}\ i\leq 0,\\
0  & \mathrm{for}\ i\leq -1,
\Ea\right.
\Eeq
where the summand $\R[-1]$ in $H^1\left(\Def(\cA ss_\infty\hspace{-0.8mm} \rar \hspace{-0.8mm}\cG raphs)\right)$ is generated by the following graph
\Beq\label{2: 3-graph}
\sum_{\sigma\in \bS_3}(-1)^\sigma
\Ba{c}
\xy
(0,-2)*{_{\sigma(2)}},
(-7,-2)*{_{\sigma(1)}},
(7,-2)*{_{\sigma(3)}},
 (-7,0)*{\circ}="a",
 (0,0)*{\circ}="b",
 (7,0)*{\circ}="c",
(0,7)*{\bu}="bu",
\ar @{-} "a";"bu" <0pt>
\ar @{-} "b";"bu" <0pt>
\ar @{-} "c";"bu" <0pt> \endxy
\Ea
\in \Def(\cA ss_\infty\hspace{-1mm}\rar\hspace{-1mm} \cG raphs).
\Eeq
and $H^0\left(\Def(\cA ss_\infty\hspace{-0.8mm} \rar \hspace{-0.8mm}\cG raphs)\right)$ is generated by
\xy
(0,2)*{_{1}},(5,2)*{_{2}},
 (0,0)*{\circ}="a",
(5,0)*{\circ}="b",
\ar @{-} "a";"b" <0pt>
\endxy.

\begin{lemma}\label{2: Lemma oh H1 of fGraphs}
$$
H^1\left(\Def(\cA ss_\infty\hspace{-0.8mm} \rar \hspace{-0.8mm}f\cG raphs)\right)=\fg\fr\ft \ \oplus \ \R[-1]
$$
\end{lemma}
\begin{proof}
As a complex $\Def(\cA ss_\infty \hspace{-0.8mm}\rar\hspace{-0.8mm}
f\cG raphs)$ is isomorphic to the tensor product of complexes
$$
\Def(\cA ss_\infty\hspace{-0.8mm} \rar \hspace{-0.8mm}\cG raphs)\ot \odot^{\bu\geq 0}(\mathsf{GC}_2[-2])
$$
so that
\Beqrn
H^1\left(\Def(\cA ss_\infty\hspace{-0.8mm} \rar \hspace{-0.8mm}f\cG raphs)\right)&=&\sum_{i\in \Z} H^i\left(\Def(\cA ss_\infty\hspace{-0.8mm} \rar \hspace{-0.8mm}\cG raphs)\right)\ot H^{-i-1}(\odot^{\bu\geq 0}\mathsf{GC}_2)\\
&=&H^1\left(\Def(\cA ss_\infty\hspace{-0.8mm} \rar \hspace{-0.8mm}f\cG raphs)\right)
\Eeqrn
as $H^{\leq -1}(\mathsf{GC}_2)=0$ and $H^{\leq -1}(\left(\Def(\cA ss_\infty\hspace{-0.8mm} \rar \hspace{-0.8mm}f\cG raphs)\right)=0$ according to Thomas Willwacher \cite{Wi,Wi2}.
\end{proof}

Note that in general the inclusion map of complexes,
$$
\Def(\cA ss_\infty\hspace{-0.9mm}\rar\hspace{-0.9mm} \cG raphs) \lon   \Def(\cA ss_\infty\hspace{-0.9mm}\rar\hspace{-0.9mm} f\cG raphs),
 $$
 induces an {\em injection}\, on cohomology,
$$
  H^i\left(\Def(\cA ss_\infty\hspace{-0.9mm}\rar\hspace{-0.9mm} \cG raphs)\right) \hook H^i\left(  \Def(\cA ss_\infty\hspace{-0.9mm}\rar\hspace{-0.9mm} f\cG raphs)\right)
$$
since $\Def(\cA ss_\infty\hspace{-0.9mm}\rar\hspace{-0.9mm} \cG raphs)$ is direct summand of $\Def(\cA ss_\infty\hspace{-0.9mm}\rar\hspace{-0.9mm} f\cG raphs)$.

\subsection{A mapping cone of the Willwacher map}
It was proven
 in \cite{Wi} that there is a degree 1 morphism of complexes,
$$
\fW : \left(\mathsf{GC}_2, \delta_\bubu\right) \lon (  \Def(\cA ss_\infty\hspace{-0.9mm}\rar\hspace{-0.9mm} f\cG raphs) ,
 \delta).
$$
which induces an {\em injection}\, on cohomology \cite{Wi, Wi2}
$$
[\fW]: H^i \left(\mathsf{GC}_2\right) \lon H^i\left(\Def(\cA ss_\infty\hspace{-0.9mm}\rar\hspace{-0.9mm} f\cG raphs)\right).
$$
The map $\fW$ sends a graph $\ga\in \mathsf{GC}_2$ (with, say, $n$ black vertices) to a
linear combination of graphs with $n$ black vertices and one white vertex,
\Beq\label{3: Willwacher map}
\fW(\ga):= \frac{1}{\# V(\ga)!}\sum_{v\in V(\ga)}
\Ba{c}\xy (0,-2)*{_{1}},
(2,4)*{_{v}},
 (0,0)*{\circ}="a",
(0,7)*{\ga}="b",
\ar @{-} "a";"b" <0pt>\endxy\Ea
=:
\Ba{c}\xy (0,-2)*{_{1}},
 (0,0)*{\circ}="a",
(0,7)*{\ga}="b",
\ar @{-} "a";"b" <0pt>\endxy\Ea
\Eeq
where $
\Ba{c}\xy (0,-2)*{_{1}},
(2,4)*{_{v}},
 (0,0)*{\circ}="a",
(0,7)*{\ga}="b",
\ar @{-} "a";"b" <0pt>\endxy\Ea$ stands for the  graph obtained by attaching
$\Ba{c}\xy
 (0,0)*{\circ}="a",
(0,7)*{}="b",
\ar @{-} "a";"b" <0pt>\endxy\Ea$ to the vertex $v$ of $\ga$; the set of edges of
$\fW(\ga)$ is ordered by putting the new edge after the edges of $\ga$.
Let $\sM\sa\sC(\fW)$ be the mapping cone of the map $\fW$, that is, the direct sum
(without the standard degree shift as the map $\fW$ has degree $+1$)
$$
\sM\sa\sC(\fW)=  \Def(\cA ss_\infty\hspace{-0.9mm}\rar\hspace{-0.9mm} f\cG raphs)  \oplus \mathsf{GC}_2
$$
equipped with the differential
$$
\Ba{rccc}
d:&  \Def(\cA ss_\infty\hspace{-0.9mm}\rar \hspace{-0.9mm}f\cG raphs^\circlearrowleft) \oplus \mathsf{GC}_2  &
\lon & \Def(\cA ss_\infty\hspace{-0.9mm}\rar\hspace{-0.9mm} f\cG raphs)  \oplus \mathsf{GC}_2 \\
& (\Ga, \ga) & \lon & (\delta\Ga + \fW(\ga), \delta_\bubu \ga).
\Ea
$$
There is a natural representation \cite{Wi} of the Lie algebra  $\mathsf{GC}_2$
on the vector space  $\Def(\cA ss_\infty\rar f\cG raphs^\circlearrowleft)$,
$$
\Ba{rccc}
\circ:&   \mathsf{GC}_2\times \Def(\cA ss_\infty\hspace{-0.9mm}\rar\hspace{-0.9mm} f\cG raphs)
& \lon & \Def(\cA ss_\infty\hspace{-0.9mm}\rar \hspace{-0.9mm}f\cG raphs) \\
& (\ga, \Ga) & \lon & \Ga \circ \ga:= \sum_{v\in V_{black}(\Ga)} \Ga\circ_v \ga
\Ea
$$
given by substitution of the graph $\ga$ into black vertices of the graph $\Ga$. This action can
 be used to make $\sM\sa\sC(\fW)$ into a Lie algebra with the brackets,
\Beq\label{3: Lie brackets in Mac(W)}
\left[(\Ga_1, \ga_1), (\Ga_2, \ga_2)\right]:= \left( [\Ga_1,\Ga_2] + \Ga_1 \circ \ga_2 -
(-1)^{|\Ga_2||\ga_1|}  \Ga_2 \circ \ga_1, \ [\ga_1,\ga_2]   \right).
\Eeq
The differential $d$ respects these brackets so that
\Beq\label{3: MapCone Lie}
\left(\sM\sa\sC(\fW), [\ ,\ ], d\right)
\Eeq
is a differential Lie algebra. For future reference we need the following

\begin{lemma}\label{3: Lemma on H1 MaC}
$H^1(\sM\sa\sC(\fW), d)=\R[-1]$.
\end{lemma}
\begin{proof} There is a short exact sequence of dg Lie algebras,
$$
0\lon  \Def(\cA ss_\infty\hspace{-0.9mm}\rar \hspace{-0.9mm}f\cG raphs) \stackrel{\al}{\lon} \sM\sa\sC(\fW)
\stackrel{\be}{\lon} \mathsf{GC}_2 \lon 0
$$
where
$$
\Ba{rccc}
\al: & \Def(\cA ss_\infty\hspace{-0.9mm}\rar \hspace{-0.9mm}\cG raphs) & \lon &
\Def(\cA ss_\infty\hspace{-0.9mm}\rar\hspace{-0.9mm} f\cG raphs)\oplus \mathsf{GC}_2\\
& \Ga & \lon & (\Ga, 0)
\Ea
$$
and
$$
\Ba{rccc}
\be: &  \Def(\cA ss_\infty\hspace{-0.9mm}\rar\hspace{-0.9mm} f\cG raphs)\oplus \mathsf{GC}_2
 & \lon & \mathsf{GC}_2\\
& (\Ga, \ga) & \lon & \ga
\Ea
$$
are the natural maps. We have, therefore, a piece of the associated  long exact sequence of
cohomology groups,
$$
H^i\hspace{-1mm}\left( \mathsf{GC}_2\right) \stackrel{[\fW]}{\rar}
 H^{i+1}\hspace{-1mm}\left( \Def(\cA ss_\infty\hspace{-1mm}\rar\hspace{-1mm} f\cG raphs) \right)
  \stackrel{[\al]}{\rar} H^{i+1}\hspace{-1mm}\left( \sM\sa\sC(\fW)\right)
\stackrel{[\al]}{\rar}  H^{i+1}\hspace{-1mm}\left( \mathsf{GC}_2\right)
\stackrel{[\fW]}{\rar} H^{i+2}\hspace{-1mm}\left( \Def(\cA ss_\infty\hspace{-1mm}\rar\hspace{-1mm} f\cG raphs)
\right)
$$
As the map $[\fW]$ is injective, we obtain
$$
 H^{i+1}\hspace{-1mm}\left( \sM\sa\sC(\fW)\right)=
 \frac{H^{i+1}\hspace{-1mm}\left( \Def(\cA ss_\infty\hspace{-1mm}\rar\hspace{-1mm} \cG raphs) \right)}
 {[\fW]( H^i\hspace{-1mm}\left( \mathsf{GC}_2\right))}.
$$
Since
$
 H^0\hspace{-1mm}\left( \mathsf{GC}_2\right)=\fg\fr\ft$ and
 $H^{1}\hspace{-1mm}\left( \Def(\cA ss_\infty\hspace{-1mm}\rar\hspace{-1mm} f\cG raphs)\right)
 =\fg\fr\ft \oplus \R[-1]$
 the claim follows.\end{proof}

%
%
%
%
%
%
%

\bip

\bip

{\large
\section{\bf Universal $nc\cG_\infty$ deformations of the standard  Gerstenhaber algebra \\ structure in $\cT_{poly}(\R^d)$}\label{Section 4: Def of ncG}
}
\subsection{Two-coloured version of $\cG ra$} Let $\cG ra=\left\{\cG ra(n)\right\}_{n\geq 1}$ be the operad defined in \S {\ref{subsection 2: Gra}}; from now one we assume that vertices of graphs from $\cG ra$ are coloured in black.
For arbitrary integers $m\geq 1$, $n\geq 0$ and $l\geq 0$ we denote by ${\sG}'_{m,n;l}$ the set
of tadpoles-free graphs\ $\{\Ga\}$ with $m$ white vertices, $n$ black
vertices and $l$ edges, such that
 \Bi
 \item[(i)] the set of white vertices, $V_\circ(\Ga)$, of $\Ga$ is equipped with a total order
 (so that in our pictures white vertices will depicted as lying on a line),
 \item[(ii)] there is a bijection  $V_\circ(\Ga) \rar [m]$ (which does not, in general, respect total orders)
 \item[(iii)] there is a bijection from the set,  $V_\bu(\Ga)$, of black vertices of $\Ga$
 to the set $[\bar{n}]=\{\bar{1},\ldots, \bar{n}\}$,
\item[(iv)] the black vertices of $\Ga$ are at least trivalent,
\item[(v)] the set of edges, $E(\Ga)$, is equipped with an orientation, i.e.\ it is totally ordered up to an even permutation.
\Ei
Note that  graphs from ${\sG}'_{m,n;l}$ can have connected components consisting of graphs with solely black vertices.
Let $\K\langle \sG_{m,{n};l}'\rangle$  be the vector space over a field $\K$ spanned by isomorphism classes,
 $[\Ga]$, of elements of $\sG_{m,{n};l}$ modulo the relation $\Ga_{opp}=-\Ga$, where
 the graph $\Ga_{opp}$ is identical to $\Ga$ except that it has the opposite orientation. Consider the following collection of $\Z$-graded $\bS$-modules,
$$
\cG ra^{\circ\bu}:=\left\{
\bigoplus_{N=m+n}\mathrm{Ind}_{\bS_N}^{\bS_m\times \bS_{\bar{n}}}\left\{\cG ra(m,{n}):=\bigoplus_{l=0}^\infty \K\langle
\sG'_{m,{n};l}\rangle[l]\right\}_{m\geq 1,n\geq 0},  \left\{\cG ra(n)\right\}_{n\geq 1}\right\}.
$$
It has a structure of a 2-coloured operad with compositions
$$
\Ba{rlcll}
\circ_i: & \cG ra(m_1,{n_1})\ot  \cG ra(m_2,{n_2}) & \lon &  \cG ra(m_1+m_2-1,{n_1+n_2}), & i\in [m_1]\\
\circ_{i}: & \cG ra(m,n_1)\ot  \cG ra(n_2) &\lon & \cG ra(m,n_1+n_2-1), & i\in [n_1]\\
\circ_{i}: & \cG ra(n_1)\ot  \cG ra(n_2) & \lon & \cG ra(n_1+n_2-1), & i\in [n_1],
\Ea
$$
given by graph substitutions as in the case of $\cG ra$.

\begin{proposition} There is a morphism of operads
$$
f: nc\cG \lon \cG ra^{\circ\bu}
$$
given on generators as follows,

\Beq\label{4: ncG to Gra2}
\xy
 <0mm,0.55mm>*{};<0mm,3.5mm>*{}**@{-},
 <0.5mm,-0.5mm>*{};<2.2mm,-2.2mm>*{}**@{-},
 <-0.48mm,-0.48mm>*{};<-2.2mm,-2.2mm>*{}**@{-},
 <0mm,0mm>*{\circ};<0mm,0mm>*{}**@{},
 <0.5mm,-0.5mm>*{};<2.7mm,-3.2mm>*{_2}**@{},
 <-0.48mm,-0.48mm>*{};<-2.7mm,-3.2mm>*{_1}**@{},
 \endxy \ \lon \
  \xy
(0,2)*{_{1}},
(5,2)*{_{2}},
 (0,0)*{\circ}="a",
(5,0)*{\circ}="b",
\endxy
,\ \ \ \ \
\Ba{c}
\xy
 <0mm,0.55mm>*{};<0mm,3.5mm>*{}**@{-},
 <0.5mm,-0.5mm>*{};<2.2mm,-2.2mm>*{}**@{-},
 <-0.48mm,-0.48mm>*{};<-2.2mm,-2.2mm>*{}**@{-},
 <0mm,0mm>*{\bu};<0mm,0mm>*{}**@{},
 <0.5mm,-0.5mm>*{};<2.7mm,-3.2mm>*{_2}**@{},
 <-0.48mm,-0.48mm>*{};<-2.7mm,-3.2mm>*{_1}**@{},
 \endxy\Ea
 \ \lon \ \xy
(0,2)*{_{1}},
(5,2)*{_{2}},
 (0,0)*{\bullet}="a",
(5,0)*{\bu}="b",
\ar @{-} "a";"b" <0pt>
\endxy
,
\ \ \ \ \ \
\Ba{c}\begin{xy}
 <0mm,0.55mm>*{};<0mm,3.5mm>*{}**@{.},
 <0.5mm,-0.5mm>*{};<2.2mm,-2.2mm>*{}**@{.},
 <-0.48mm,-0.48mm>*{};<-2.2mm,-2.2mm>*{}**@{-},
 <0mm,0mm>*{\bu};<0mm,0mm>*{}**@{},
 <0.5mm,-0.5mm>*{};<2.8mm,-3.3mm>*{_2}**@{},
 <-0.48mm,-0.48mm>*{};<-2.7mm,-3.3mm>*{_1}**@{},
 \end{xy}\Ea\  \lon \
 \Ba{c}\xy (0,6)*{_1},(0,-2)*{_{2}},
 (0,4)*{\bu}="a", (0,0)*{\circ}="b",
\ar @{-} "a";"b" <0pt>
\endxy\Ea
\Eeq
\end{proposition}
\begin{proof} We have to check that the map $f$ respects relations (\ref{2: assoc relations}),
(\ref{2: Jacobi relations}) and (\ref{2: noncomm Gerstenhaber condition}). For example,
\Beqrn
f\left(
\Ba{c}
\resizebox{0.9cm}{!}{
\begin{xy}
 <0mm,0mm>*{\bu};<0mm,0mm>*{}**@{-},
 <0mm,0.49mm>*{};<0mm,3.5mm>*{}**@{.},
 <0.49mm,-0.49mm>*{};<1.9mm,-1.9mm>*{}**@{.},
 <-0.5mm,-0.5mm>*{};<-1.9mm,-1.9mm>*{}**@{-},
 <2.3mm,-2.3mm>*{\circ};<-2.3mm,-2.3mm>*{}**@{},
 <1.8mm,-2.8mm>*{};<0mm,-4.9mm>*{}**@{.},
 <2.8mm,-2.9mm>*{};<4.6mm,-4.9mm>*{}**@{.},
   <0.49mm,-0.49mm>*{};<-2.8mm,-3.2mm>*{_1}**@{},
   <-1.8mm,-2.8mm>*{};<0.1mm,-6.2mm>*{_2}**@{},
   <-2.8mm,-2.9mm>*{};<5.3mm,-6.2mm>*{_3}**@{},
 \end{xy}}
 \Ea
 -
\Ba{c}
\resizebox{0.8cm}{!}{
\begin{xy}
 <0mm,0mm>*{\circ};<0mm,0mm>*{}**@{},
 <0mm,0.49mm>*{};<0mm,3.5mm>*{}**@{.},
 <0.49mm,-0.49mm>*{};<1.9mm,-1.9mm>*{}**@{.},
 <-0.5mm,-0.5mm>*{};<-1.9mm,-1.9mm>*{}**@{.},
 <-2.3mm,-2.3mm>*{\bu};<-2.3mm,-2.3mm>*{}**@{},
 <-1.8mm,-2.8mm>*{};<0mm,-4.9mm>*{}**@{.},
 <-2.8mm,-2.9mm>*{};<-4.6mm,-4.9mm>*{}**@{-},
   <0.49mm,-0.49mm>*{};<2.7mm,-3.2mm>*{_3}**@{},
   <-1.8mm,-2.8mm>*{};<0.4mm,-6.2mm>*{_2}**@{},
   <-2.8mm,-2.9mm>*{};<-5.1mm,-6.2mm>*{_1}**@{},
 \end{xy}}\Ea
 -
\Ba{c}
\resizebox{0.8cm}{!}{
\begin{xy}
 <0mm,0mm>*{\circ};<0mm,0mm>*{}**@{},
 <0mm,0.49mm>*{};<0mm,3.5mm>*{}**@{.},
 <0.49mm,-0.49mm>*{};<1.9mm,-1.9mm>*{}**@{.},
 <-0.5mm,-0.5mm>*{};<-1.9mm,-1.9mm>*{}**@{.},
 <2.3mm,-2.3mm>*{\bu};<-2.3mm,-2.3mm>*{}**@{},
 <1.8mm,-2.8mm>*{};<0mm,-4.9mm>*{}**@{-},
 <2.8mm,-2.9mm>*{};<4.6mm,-4.9mm>*{}**@{.},
   <0.49mm,-0.49mm>*{};<-2.8mm,-3.2mm>*{_2}**@{},
   <-1.8mm,-2.8mm>*{};<0.1mm,-6.2mm>*{_1}**@{},
   <-2.8mm,-2.9mm>*{};<5.3mm,-6.2mm>*{_3}**@{},
 \end{xy}}
 \Ea
 \right) &=&
\Ba{c}\xy (0,6)*{_1},(0,-2)*{_{2}},
 (0,4)*{\bu}="a", (0,0)*{\circ}="b",
\ar @{-} "a";"b" <0pt>
\endxy\Ea
\
 \circ_2 \
\Ba{c} \xy
(0,2)*{_{1}},
(3,2)*{_{2}},
 (0,0)*{\circ}="a",
(3,0)*{\circ}="b",
\endxy\Ea \ - \ \Ba{c} \xy
(0,2)*{_{1}},
(3,2)*{_{2}},
 (0,0)*{\circ}="a",
(3,0)*{\circ}="b",
\endxy\Ea
\circ_1
\Ba{c}\xy (0,6)*{_1},(0,-2)*{_{2}},
 (0,4)*{\bu}="a", (0,0)*{\circ}="b",
\ar @{-} "a";"b" <0pt>
\endxy\Ea
\ - \ \Ba{c} \xy
(0,2)*{_{1}},
(3,2)*{_{2}},
 (0,0)*{\circ}="a",
(3,0)*{\circ}="b",
\endxy\Ea
\circ_2
\Ba{c}\xy (0,6)*{_1},(0,-2)*{_{2}},
 (0,4)*{\bu}="a", (0,0)*{\circ}="b",
\ar @{-} "a";"b" <0pt>
\endxy\Ea
\\
&=&
\Ba{c}\xy (0,6)*{_1},(0,-2)*{_{2}}, (3,-2)*{_{3}},
 (0,4)*{\bu}="a", (0,0)*{\circ}="b",  (3,0)*{\circ},
\ar @{-} "a";"b" <0pt>
\endxy\Ea
\ + \
\Ba{c}\xy (0,6)*{_1},(0,-2)*{_{3}}, (-3,-2)*{_{2}},
 (0,4)*{\bu}="a", (0,0)*{\circ}="b",  (-3,0)*{\circ},
\ar @{-} "a";"b" <0pt>
\endxy\Ea
\ - \
\Ba{c}\xy (0,6)*{_1},(0,-2)*{_{2}}, (3,-2)*{_{3}},
 (0,4)*{\bu}="a", (0,0)*{\circ}="b",  (3,0)*{\circ},
\ar @{-} "a";"b" <0pt>
\endxy\Ea
\ - \
\Ba{c}\xy (0,6)*{_1},(0,-2)*{_{3}}, (-3,-2)*{_{2}},
 (0,4)*{\bu}="a", (0,0)*{\circ}="b",  (-3,0)*{\circ},
\ar @{-} "a";"b" <0pt>
\endxy\Ea\\
&=&0.
\Eeqrn
Analogously one checks all other relations. \end{proof}

\begin{theorem} The deformation complex, $\Def(nc\cG_\infty\rar \cG ra^{\circ\bu})$, of the morphism
$$
f_o: nc\cG_\infty \xrightarrow{proj} nc\cG \xrightarrow{ \ f\  }  \cG ra^{\circ\bu}
$$
is isomorphic as a dg Lie algebra to $\mathsf{MaC}(\fW)$.
\end{theorem}
\begin{proof} As a graded vector space $\Def(nc\cG_\infty\rar \cG ra^{\circ\bu})$
is identical to the space of homomorphisms, $\Hom_\bS(E,\cG ra^{\circ\bu})[-1]$, of $\bS$-modules,
where $E=\{E(N)\}$ is the $\bS$-module of generators of $nc\cG_\infty$. The latter $\bS$-module
splits as a direct sum,
$$
E(N)=E_1(N) \oplus E_2(N),
$$
where $E_1(N)$ is spanned as a vector space by corollas (\ref{2: Lie_inf corolla}) and hence is given by
$$
E_1(N)=sgn_N[2n-3]
$$
where $sgn_N$ is the the one-dimensional sign representation of $\bS_N$. The $\bS_N$-module
 $E_2(N)$ is spanned by corollas (\ref{2: OCHA corolla}) and hence equals
 $$
 E_2(N)=\bigoplus_{N=m+n\atop m\geq 1, n\geq 0} \mathrm{Ind}_{\bS_N}^{\bS_m\times \bS_n} \K[\bS_m]\ot sgn_n[2n+m-2].
 $$
Therefore, we have an isomorphism of graded vector spaces
\Beqrn
\Def(nc\cG_\infty\rar \cG ra^{\circ\bu})&=&\bigoplus_{N}\Hom_\bS\left(E_2(N),\cG ra^{\circ\bu}(N)\right)[-1] \oplus \bigoplus_{N}\Hom_\bS\left(E_1(N),\cG ra^{\circ\bu}(N)\right)[-1]\\
&=&  \Def(\cA ss_\infty\hspace{-0.9mm}\rar\hspace{-0.9mm} f\cG raphs)  \oplus \mathsf{GC}_2\\
&=& \mathsf{MaC}(\fW)
\Eeqrn
One reads the Lie algebra structure in  $\Def(nc\cG_\infty\rar \cG ra^{\circ\bu})$ from the differential (\ref{2: d on Lie corollas}) and (\ref{2: differential on Konts corollas})  and easily concludes that
it is given precisely by the Lie bracket $[\ ,\ ]$ given in (\ref{3: Lie brackets in Mac(W)}).
Next, there is a 1-1 correspondence between Maurer-Cartan elements, $\Ga$,
$$
[\Ga, \Ga]=0,
$$
and morphisms of operads $nc\cG_\infty \rar  \cG ra^{\circ\bu}$ (cf.\ \cite{MV}). The morphism $f_0$ is represented  by the following Maurer-Cartan element,
\Beq\label{4: Ga_0}
\Ga_0=  \left( \xy
 (0,0)*{\circ}="a",
(3,0)*{\circ}="b",
\endxy\ +\ \xy
 (0,4)*{\bu}="a", (0,0)*{\circ}="b",
\ar @{-} "a";"b" <0pt>
\endxy\ \ , \ \  \xy
 (0,0)*{\bullet}="a",
(4,0)*{\bu}="b",
\ar @{-} "a";"b" <0pt>
\endxy\right)
\Eeq
so that the differential in  $\Def(nc\cG_\infty\rar \cG ra^{\circ\bu})$
is given by $[\Ga_0,\ ]$ and hence coincides precisely with the differential $d$ in
$\mathsf{MaC}(\fW)$. The theorem is proven.
\end{proof}

\subsection{A canonical representations of $\cG ra^{\circ\bu}$  in polyvector fields and an exotic $nc\cG_\infty$ structure}
There is a representation of the two-colored  operad $\cG ra^{\circ\bu}$
in the two-coloured endomorphism operad, $\cE nd_{\cT_{poly}(R^d), \cT_{poly}(R^d)}$, of two copies of the space  $\cT_{poly}(R^d)$ given by formulae which are
completely analogous to (\ref{3: Gra representation}). Hence there is an induced of morphism
of dg Lie algebras
$$
\mathsf{MaC}(\fW)=\Def(nc\cG_\infty\rar \cG ra^{\circ\bu}) \lon \Def(nc\cG_\infty\rar \cE nd_{\cT_{poly}(R^d), \cT_{poly}(R^d)}).
$$
The dg Lie algebra $\Def(nc\cG_\infty\rar \cE nv_{\cT_{poly}(R^d), \cT_{poly}(R^d)})$ describes
$nc\cG_\infty$ deformation of the standard Gerstenhaber algebra structure on $ \cT_{poly}(R^d)$.
The dg Lie algebra $\mathsf{MaC}(\fW)$ controls, therefore, {\em universal}\, deformations of this structure, i.e\ the ones which make sense in any dimension $d$.

\sip

In particular any Maurer-Cartan element,
$$
d\Ga + \frac{1}{2}[\Ga,\Ga]=0
$$
in the dg Lie algebra $\mathsf{MaC}(\fW)$ gives us a universal $nc\cG_\infty$-structure in  $\cT_{poly}(R^d)$.
Such a structure can be viewed as a deformation of the standard Gerstenhaber algebra structure
(corresponding to the graph (\ref{4: Ga_0})) as the above equation can be rewritten as
$$
[\Ga_0+\Ga,\Ga_0+\Ga]=0.
$$
The  dg Lie algebra   $\mathsf{MaC}(\fW)$ is naturally filtered by the number of black and white vertices. We assume from now on that  $\mathsf{MaC}(\fW)$ is completed with respect to this filtration. Then there is a well-defined action of degree zero elements, $g$, of  $\mathsf{MaC}(\fW)$ on the set of Maurer-Cartan elements,
$$
    \mathsf{\Ga} \lon \Ga^g:= e^{ad_g}  \mathsf{\Ga} - \frac{ e^{ad_g} - 1}{ad_g}d g.
$$
The orbits of this action are $nc\cG_\infty$-isomorphism classes of universal $nc\cG_\infty$ structures on polyvector fields.

\sip

Infinitesimal homotopy non-trivial $nc\cG_\infty$ deformations of the standard Gerstenhaber algebra structure on polyvector fields are classified by the cohomology group
$H^1(\mathsf{MaC}(\fW))$. Lemma~{\ref{3: Lemma on H1 MaC}} says that there exists at most one homotopy non-trivial universal  $nc\cG_\infty$ deformation of the standard Gerstenhaber algebra structure on polyvector fields. The associated Maurer-Cartan element in $\mathsf{MaC}(\fW)$ was
given explicitly in \cite{A} in term of periods over the compactified configuration spaces
$\overline{C}_{\bu,\bu}(\C)$,
\Beq\label{4: Exotic ncG}
\Ga_0+\Ga=\left(\sum_{m\geq 1, n\geq 0}\ \ \ \sum_{\Ga\in {\sG}''_{m,n;2n+m-2}} \underset{\overline{C}_{m,n}(\C)}{\int}\hspace{-1mm}\Omega_\Ga\  \Ga,\ \
\ \  \xy
 (0,0)*{\bullet}="a",
(4,0)*{\bu}="b",
\ar @{-} "a";"b" <0pt>
\endxy   \right)
\Eeq
where
\Bi
\item ${\sG}''_{m,n;2n+m-2}$ is the set of equivalence classes of graphs from $\sG'_{m,n;2n+m-2}$ which are linearly independent in the space $\K\langle{\sG}'_{m,n;2n+m-2}\rangle$ and have no tadpoles;

\item
 $\Omega_\Ga:=\bigwedge_{e\in Edges(\Ga)}{\pi}_e^*(\om)$ ,

\item  for an edge
$e\in Edges(\Ga)$
beginning at a vertex (of any colour) labelled by $i$ and ending at a vertex (of any colour) labelled by $j$, $p_e$ is the natural surjection
$$
\Ba{rccc}
\pi_{e}: & C_{n,m}(\C) & \lon & C_2(\C)=S^1\\
& (z_1,\ldots, z_i, \ldots, z_j,\ldots, z_{n+m}) & \lon & \frac{z_i - z_j}{|z_i-z_j|}.
\Ea
$$
\item
The $1$-form $\om:=\frac{1}{2\pi}dArg(z_i-z_j)$ is the standard homogenous volume form
on $S^1$ normalized so that $\int_{S^1}\om=1$.
\Ei
The lowest (in total number of vertices) term
in $\Ga$  is given by the graph (\ref{2: 3-graph}) whose weight is equal to $1/3$. Hence
Lemma~{\ref{3: Lemma on H1 MaC}} and \cite{A} imply the following

\begin{theorem}\label{4: Theorem on 2 ncG-infty}
Up to $nc\cG_\infty$ isomorphisms, there are only two different universal  $nc\cG_\infty$ structures on polyvector fields,
the standard Gerstenhaber one corresponding to the Maurer-Cartan element (\ref{4: Ga_0})
and the exotic one given by (\ref{4: Exotic ncG}).
\end{theorem}

\bip

\bip


{\large
\section{\bf No-Go Theorem}\label{Section 6: No-Go}
}

\bip

\subsection{A class of universal  $\cA ss_\infty$ structures on Poisson manifolds} For any degree $2$ element $\hbar \pi$ in $\cT_{poly}(\R^d)[[\hbar]]$ the operad
$f\cG raphs$ admits a canonical representation
$$
\rho^\pi: f\cG raphs \lon \cE nd^{cont}_{\cT_{poly}(\R^d)[[\hbar]]}
$$
which sends a graph $\Ga$ from  $\cG raphs$ with, say, $m$ white vertices and $n$ black vertices
into a continuous (in the $\hbar$-adic topology) operator $\rho(\Ga)\in \Hom (\ot^m\cT_{poly}(\R^d), \cT_{poly}(\R^d)[[[\hbar]]$ which is constructed exactly as in the formula (\ref{3: Gra representation}) except that black vertices are decorated by the element $\hbar\pi$. (From now on we take our operad $f\cG raphs$
to be completed with respect to the filtration by the number of black vertices; hence we need to use a degree zero formal parameter $\hbar$ to ensure convergence of operators $\rho(\Ga)$ in the $\hbar$-adic topology.)

\sip

The Lie algebra $\mathsf{GC}_2$ acts (on the right) on the operad $\cG raphs$,
$$
\Ba{rccc}
R: &  f\cG raphs \times \mathsf{GC}_2 & \lon &  \cG raphs\\
   &     (\Ga, \ga) & \lon & \Ga \cdot \ga
\Ea
$$
by operadic derivations, where $\Ga\cdot  \ga$ is obtained from $\Ga$ by inserting $\ga$ into black vertices
\cite{Wi}. Let $I'_\bubu$ be the operadic ideal in  $\cG raphs$ generated by graphs of the form
$\Ga\cdot \bubu$. There is natural projection map of operads,
$$
 f\cG raphs \lon  f\cG raphs':= f\cG raphs/I'_\bubu,
$$
and, for $\pi$ being a (graded) Poisson structure on $\R^d$, that is, for $\pi$ satisfying $[\pi,\pi]_{S}=0$, the canonical representation $\rho^\pi$ factors through this projection,
$$
\rho^\pi: f\cG raphs \lon  f\cG raphs' \lon \cE nd_{\cT_{poly}(\R^d)}.
$$
The induced representation  $f\cG raphs' \rar \cE nd_{\cT_{poly}(\R^d)}$ we denote by the same letter $\rho^\pi$. It induces in turn a map of Lie algebras,
$$
\Def(\cA ss_\infty\hspace{-0.9mm}\rar \hspace{-0.9mm}f\cG raphs)\lon \Def(\cA ss_\infty\hspace{-0.9mm}\rar \hspace{-0.9mm}f\cG raphs')\lon \Def(\cA ss_\infty\hspace{-0.9mm}\rar \hspace{-0.9mm} \cE nd_{\cT_{poly}(\R^d)})= \mathsf{CoDer}(\ot^{\bu \geq 1}(\cT_{poly}(\R^d)[1])
$$
where
$$
\mathsf{CoDer}(\ot^{\bu \geq 1}(\cT_{poly}(\R^d)[1])=\bigoplus_{m\geq 1}
\mathsf{Hom}(\ot^m \cT_{poly}(\R^d),
\cT_{poly}(\R^d))[1-m]
$$
is the Gerstenhaber Lie algebra of coderivations of the tensor
coalgebra $\ot^{\bu\geq 1}(\cT_{poly}(\R^d)[1])$. Hence any Maurer-Cartan element $\Ga$,
$$
[\Ga,\Ga]=0,
$$
 in the Lie algebra  $\Def(\cA ss_\infty\hspace{-0.9mm}\rar \hspace{-0.9mm}f\cG raphs')$ induces,
 for any fixed Poisson structure on $\R^d$, a universal $\cA ss_\infty$ algebra structure on $\cT_{poly}(\R^d)$. Moreover, two such universal $\cA ss_\infty$ structures, $\Ga_1$ and $\Ga_2$, are universally $\cA ss_\infty$ isomorphic if and only if there exists a degree zero element $h\in
 \Def(\cA ss_\infty\hspace{-0.9mm}\rar \hspace{-0.9mm}f\cG raphs')$ such that
 $$
 \Ga_2=e^{\mathrm{ad}_h}\Ga_1,
 $$
 where $\mathrm{ad}_h$ stands for the adjoint action. Note that, due to the filtrations of the Lie algebra
 $\Def(\cA ss_\infty\hspace{-0.9mm}\rar \hspace{-0.9mm}f\cG raphs')$ by the numbers of white and black vertices, there is no convergence problem in taking the exponent of  $\mathrm{ad}_h$.

 \sip

 It is easy to see that
 $$
 \Ga_0= \wiwi + \wibu
 $$
is a Maurer-Cartan element in $\Def(\cA ss_\infty\hspace{-0.9mm}\rar \hspace{-0.9mm}f\cG raphs')$
(but {\em not}\, in $\Def(\cA ss_\infty\hspace{-0.9mm}\rar \hspace{-0.9mm}f\cG raphs)!)$, and the associated
$\cA ss_\infty$ structure in $(\cT_{poly}(\R^d), \pi)$ is the standard structure of a Poisson complex,
that is, a wedge product $\wedge$ (corresponding to the graph $\wiwi$) and the differential $d_\pi=[\hbar\pi,\ ]_{\mathrm{S}}$ (corresponding to the graph $\wibu$). Hence $\Ga_0$ makes $\Def(\cA ss_\infty\hspace{-0.9mm}\rar \hspace{-0.9mm}f\cG raphs')$ into a complex with the differential
$d:=[\Ga_0,\ ]$. It is clear that the natural projection of Lie algebras
$$
p: \Def(\cA ss_\infty\hspace{-0.9mm}\rar \hspace{-0.9mm}f\cG raphs)\lon \Def(\cA ss_\infty\hspace{-0.9mm}\rar \hspace{-0.9mm}f\cG raphs')
$$
is compatible with the differentials.

\sip

Let $\ga$ be a degree zero cycle in the graph complex $\mathsf{GC}_2$, representing some cohomology class
from $\fg\fr\ft$. Then
$$
\Ga_0^\ga:=\Ga_0\cdot e^\ga=  \wiwi + \wibu + \Ba{c}\xy (0,6)*{^{\ga}},
 (0,4)*{\bu}="a", (0,0)*{\circ}="b",
\ar @{-} "a";"b" <0pt>\endxy\Ea   +\frac{1}{2!}\Ba{c}\xy (0,6)*{^{\ga\cdot \ga}},
 (0,4)*{\bu}="a", (0,0)*{\circ}="b",
\ar @{-} "a";"b" <0pt>
\endxy\Ea  +\ldots
$$
is again a Maurer-Cartan element in $\Def(\cA ss_\infty\hspace{-0.9mm}\rar \hspace{-0.9mm}f\cG raphs')$.
The associated
$\cA ss_\infty$ structure in $(\cT_{poly}(\R^d), \pi)$ consists of the standard
wedge product $\wedge$ (corresponding to the graph $\wiwi$) and the differential $d_{g(\pi)}=[g(\hbar\pi),\ ]_{\mathrm{S}}$, where $g=\exp(\ga)$ is the element of the group $GRT$
corresponding to $\ga$. The element in the difference $\Ga_0^\ga -\Ga_0$ with lowest number of vertices is
$$
 \Ba{c}\xy (0,6)*{^{\ga}},
 (0,4)*{\bu}="a", (0,0)*{\circ}="b",
\ar @{-} "a";"b" <0pt>\endxy\Ea
$$
It defines a cycle in both complexes $\Def(\cA ss_\infty\hspace{-0.9mm}\rar \hspace{-0.9mm}f\cG raphs)$
and $\Def(\cA ss_\infty\hspace{-0.9mm}\rar \hspace{-0.9mm}f\cG raphs')$. It is shown in \cite{Wi} that
$\Ba{c}\xy (0,6)*{^{\ga}},
 (0,4)*{\bu}="a", (0,0)*{\circ}="b",
\ar @{-} "a";"b" <0pt>\endxy\Ea$ is {\em not}\, a coboundary in $\Def(\cA ss_\infty\hspace{-0.9mm}\rar \hspace{-0.9mm}f\cG raphs)$, it therefore defines a non-trivial cohomology class in
$H^1(\Def(\cA ss_\infty\hspace{-0.9mm}\rar \hspace{-0.9mm}f\cG raphs))$.

\subsubsection{\bf Lemma}\label{5: Lemma on injection of 1-cohomology}
{\em  For any $[\ga]\in \fg\fr\ft$ an associated cycle $\Ba{c}\xy (0,6)*{^{\ga}},
 (0,4)*{\bu}="a", (0,0)*{\circ}="b",
\ar @{-} "a";"b" <0pt>\endxy\Ea$ is\, {\em not}\, a coboundary in $\Def(\cA ss_\infty\hspace{-0.9mm}\rar \hspace{-0.9mm}f\cG raphs')$, that is, it defines a non-trivial cohomology class in
$H^1(\Def(\cA ss_\infty\hspace{-0.9mm}\rar \hspace{-0.9mm}f\cG raphs'))$. In fact the natural map
$$
[p]: H^1(\Def(\cA ss_\infty\hspace{-0.9mm}\rar \hspace{-0.9mm}f\cG raphs)) \lon
H^1(\Def(\cA ss_\infty\hspace{-0.9mm}\rar \hspace{-0.9mm}f\cG raphs'))
$$
is an injection.}

\bip

Let us first prove the following corollary to this lemma, and then the lemma itself.

\subsubsection{{\bf No-go theorem}} {\em For any $[\ga]\in \fg\fr\ft$, the Maurer-Cartan elements\, $\Ga_0$ and $\Ga_0^\ga$ are not gauge equivalent in the Lie algebra $\Def(\cA ss_\infty\hspace{-0.9mm}\rar \hspace{-0.9mm}f\cG raphs')$. Equivalently, the universal $\cA ss_\infty$ structures in $\cT_{poly}(\R^d)[[\hbar]]$  corresponding to these elements are not universally $\cA ss_\infty$ isomorphic.}

\begin{proof} Comparing the terms with the same number of black and white vertices in the equation
$$
\Ga_0^\ga=e^{\mathrm{ad}_h}\Ga_0,
$$
we immediately see that
$$\Ba{c}\xy (0,6)*{^{\ga}},
 (0,4)*{\bu}="a", (0,0)*{\circ}="b",
\ar @{-} "a";"b" <0pt>\endxy\Ea=[h',  \wiwi + \wibu] = -dh'
$$
for some summand $h'$ in $h$. This contradicts Lemma~{\ref{5: Lemma on injection of 1-cohomology}}.
\end{proof}

To prove Lemma~{\ref{5: Lemma on injection of 1-cohomology}} we need the following

\subsubsection{\bf Lemma}\label{5: Lemma on reprent of grt witn no black vertices}{\em For any $[\ga]\in \fg\fr\ft$, an associated cycle  $\Ba{c}\xy (0,6)*{^{\ga}},
 (0,4)*{\bu}="a", (0,0)*{\circ}="b",
\ar @{-} "a";"b" <0pt>\endxy\Ea$ in the complex $\Def(\cA ss_\infty\hspace{-0.9mm}\rar \hspace{-0.9mm}f\cG raphs)$ is cohomologous to an element $\ga_{w}$ which has no black vertices.}

\begin{proof}
Let us represent the total differential $\delta$ in $\Def(\cA ss_\infty\hspace{-0.9mm}\rar \hspace{-0.9mm}f\cG raphs)$ as a sum of two differentials (see (\ref{3: differential in Def(Ass--Graphs)}))
$$
\delta=\delta_\wiwi + \delta'.
$$
The cohomology of the complex $(\Def(\cA ss_\infty\hspace{-0.9mm}\rar \hspace{-0.9mm}\cG raphs), \delta')$
(which contains elements of the form $\Ba{c}\xy (0,6)*{^{\ga}},
 (0,4)*{\bu}="a", (0,0)*{\circ}="b",
\ar @{-} "a";"b" <0pt>\endxy\Ea$ and  is a {\em direct}\, summand of the full complex $(\Def(\cA ss_\infty\hspace{-0.9mm}\rar \hspace{-0.9mm}f\cG raphs), \delta')$     was computed in \cite{LV} (see also
Proposition 5 in \cite{Wi}). We need from that computation only the following fact: any $\delta'$-cocycle
in  $\Def(\cA ss_\infty\hspace{-0.9mm}\rar \hspace{-0.9mm}\cG raphs)$ which contains at least one black vertex is $\delta'$-exact. As
$$
\delta \Ba{c}\xy (0,6)*{^{\ga}},
 (0,4)*{\bu}="a", (0,0)*{\circ}="b",
\ar @{-} "a";"b" <0pt>\endxy\Ea =\delta_\wiwi \Ba{c}\xy (0,6)*{^{\ga}},
 (0,4)*{\bu}="a", (0,0)*{\circ}="b",
\ar @{-} "a";"b" <0pt>\endxy\Ea=0,
$$
we conclude that $\delta'\Ba{c}\xy (0,6)*{^{\ga}},
 (0,4)*{\bu}="a", (0,0)*{\circ}="b",
\ar @{-} "a";"b" <0pt>\endxy\Ea=0$ and hence $\Ba{c}\xy (0,6)*{^{\ga}},
 (0,4)*{\bu}="a", (0,0)*{\circ}="b",
\ar @{-} "a";"b" <0pt>\endxy\Ea=-\delta' \ga_\circ$ for some degree zero graph $\ga_\circ$ in
$\Def(\cA ss_\infty\hspace{-0.9mm}\rar \hspace{-0.9mm}\cG raphs)$; in fact, it is easy to see that
$\ga_\circ$ is $\ga$ with every black vertex labelled by, say, $1$ made white (remember that $\ga$ is symmetrized over numerical labellings of vertices so that nothing depends on the choice of a particular label in this construction of $\ga_\circ$).  We can, therefore, write,
$$
\Ba{c}\xy (0,6)*{^{\ga}},
 (0,4)*{\bu}="a", (0,0)*{\circ}="b",
\ar @{-} "a";"b" <0pt>\endxy\Ea= -(\delta_\wiwi + \delta') \ga_\circ +\delta_\wiwi \ga_\circ.
$$
If $\delta_\wiwi \ga_\circ$ contains black vertices, then again
$$
\delta\, \delta_\wiwi \ga_\circ =\delta'(\delta_\wiwi \ga_\circ)=0\ \ \ \Rightarrow \ \ \ \delta_\wiwi \ga_\circ=-\delta' \ga_{\wiwi}
$$
and hence
$$
\Ba{c}\xy (0,6)*{^{\ga}},
 (0,4)*{\bu}="a", (0,0)*{\circ}="b",
\ar @{-} "a";"b" <0pt>\endxy\Ea= -\delta( \ga_0 + \ga_\wiwi) + \delta_\wiwi\ga_\wiwi.
$$
Continuing this process we finally obtain an equality
\Beq\label{5: ga_w}
\Ba{c}\xy (0,6)*{^{\ga}},
 (0,4)*{\bu}="a", (0,0)*{\circ}="b",
\ar @{-} "a";"b" <0pt>\endxy\Ea= -\delta( \ga_0 + \ga_\wiwi + \ldots + \ga^{max}_{\circ\ldots \circ}) + \delta_\wiwi\ga_{\circ\ldots \circ}
\Eeq
where $\ga_w:= \delta_\wiwi\ga^{max}_{\circ\ldots \circ}$ has no black vertices.
\end{proof}

\bip

{\em Proof of Lemma}\, {\ref{5: Lemma on injection of 1-cohomology}}. Since
$$
H^1(\Def(\cA ss_\infty\hspace{-0.9mm}\rar \hspace{-0.9mm}f\cG raphs))=
 H^1(\Def(\cA ss_\infty\hspace{-0.9mm}\rar \hspace{-0.9mm}\cG raphs))
 $$
and since $\Def(\cA ss_\infty\hspace{-0.9mm}\rar \hspace{-0.9mm}\cG raphs')$ is a direct summand
of $\Def(\cA ss_\infty\hspace{-0.9mm}\rar \hspace{-0.9mm}f\cG raphs')$, it is enough to study the natural projection map
$$
p: \Def(\cA ss_\infty\hspace{-0.9mm}\rar \hspace{-0.9mm}\cG raphs) \lon \Def(\cA ss_\infty\hspace{-0.9mm}\rar \hspace{-0.9mm}\cG raphs').
$$
Consider the following {\em direct}\, summands, $C$ and $C'$, of both complexes of the form
$$
\Ker d \bigcap \{\mbox{Subspace of graphs with no black vertices}\}
$$
As the ideal used to construct the quotient operad $\cG raphs'$ out of $Graphs$ consists of graphs with at least two black vertices, we conclude that the map $p$ sends
$C$ isomorphically to $C'$. Then Lemma {\ref{5: Lemma on reprent of grt witn no black vertices}}
(and its obvious analogue for the graph (\ref{2: 3-graph})) implies the required result. \hfill $\Box$

\subsection{Quotient mapping cone}\label{5: quotient maping cone} Let $I_\bubu$ be the ideal in the operad $\cG ra$ generated
by the graph
 \xy
(0,2)*{_{1}},
(5,2)*{_{2}},
 (0,0)*{\bullet}="a",
(5,0)*{\bu}="b",
\ar @{-} "a";"b" <0pt>
\endxy, let $\cG ra':=\cG ra/I_\bubu$, and let
$$
\mathsf{GC}_2':=\Def\left(\caL ie\{1\}_\infty \stackrel{0}{\rar} \cG ra'\right)
$$
be the deformation complex of the zero map (this is just a Lie algebra). There is an induced Willwacher map
$$
\fW': \mathsf{GC}_2' \lon  \Def(\cA ss_\infty\hspace{-0.9mm}\rar \hspace{-0.9mm}\cG raphs')
$$
and hence an associated Lie algebra structure on the {\em quotient}\, mapping cone,
$$
\mathsf{MaC}(\fW'):= \Def(\cA ss_\infty\hspace{-0.9mm}\rar \hspace{-0.9mm}\cG raphs')\ \oplus\
\mathsf{GC}_2'.
$$
There is a natural surjection of Lie algebras,
\Beq\label{5: surj of mapping cones}
S: \mathsf{MaC}(\fW)  \lon \mathsf{MaC}(\fW').
\Eeq
For future reference we make an evident observation that our class of universal $\cA ss_\infty$ structures on polyvector fields can be identified with a class of Maurer-Cartan elements
of the quotient mapping cone $\mathsf{MaC}(\fW')$ which have the form $(\Ga,0)$
for some $\Ga\in \Def(\cA ss_\infty\hspace{-0.9mm}\rar \hspace{-0.9mm}\cG raphs')$.

{\large
\section{\bf Proof of the main theorem}\label{Section 5: proof of the main theorem}
}

\subsection{$nc\cG_\infty$ isomorphisms of $nc\cG_\infty$ algebras} As the two-coloured operad
$\cG ra^{\circ\bu}$ has a canonical representation in the space of polyvector fields $\cT_{poly}(\R^d)$,
any morphism of operads
$$
F: nc\cG_\infty \lon \cG ra^{\circ\bu}
$$
induces a universal $nc\cG_\infty$ structure in $\cT_{poly}(\R^d)$. On the other hand, we proved in the previous section that there is a one-to-one correspondence between such morphisms
$F$  and degree 1 elements,
 $$
 \mathsf{\Ga} =(\Ga, \ga)
 $$
 in the Lie algebra $\mathsf{MaC}(\fW)= \Def(\cA ss_\infty\rar f\cG raphs) \oplus \mathsf{GC}_2$
 satisfying the Maurer-Cartan condition
 $$
 [ \mathsf{\Ga},  \mathsf{\Ga}]=\left( [\Ga,\Ga] + 2 \Ga \circ \ga, \ [\ga,\ga]   \right)=0.
 $$
Two universal $nc\cG_\infty$ structures corresponding to Maurer-Cartan elements $\mathsf{\Ga}$ and $\mathsf{\Ga}'$
are $nc\cG_\infty$-isomorphic if and only if the Maurer-Cartan elements $\mathsf{\Ga}$ and $\mathsf{\Ga}'$ are gauge equivalent, that is,
\Beq\label{6: Gauge transf in MaC}
\mathsf{\Ga}'=e^{\mathrm{ad}_\mathsf{H}}\mathsf{\Ga}
\Eeq
for some degree zero element $\mathsf{H}=(H,h)$ in $\mathsf{MaC}(\fW)$.

\subsection{$nc\cG_\infty$ structures versus $\cA ss_\infty$ structures on (affine) Poisson manifolds} Let us denote by $\cM\cC$
the set of all Maurer-Cartan elements in the Lie algebra $\mathsf{MaC}(\fW)$.
By Theorem {\ref{4: Theorem on 2 ncG-infty}}, any  element $\Ga\in \cM\cC$ is gauge equivalent either to (\ref{4: Ga_0})
or to (\ref{4: Exotic ncG}). Both these Maurer-Cartan elements belong to the subset
$\cM\cC_{\cA ss}\subset \cM\cC$ consisting of elements of the form $(\Ga, \bubu)$ for some
$\Ga\in \Def(\cA ss_\infty\rar f\cG raphs)$. As projection (\ref{5: surj of mapping cones}) sends such and element into a Maurer-Cartan element in $\mathsf{MaC}(\fW')$ of the form $(\Ga,0)$, the subset $\cM\cC_{\cA ss}\subset \cM\cC$ gives us universal $\cA ss_\infty$ structures on polyvector fields. We are interested now in the gauge transformations of the set $\cM\cC$ which preserve the subset $\cM\cC_{\cA ss}$, as such transformations can sometimes induce  (via the surjection (\ref{5: surj of mapping cones})) $\cA ss_\infty$ isomorphisms of our class of universal $\cA ss_\infty$ structures on polyvector fields. It is clear that the gauge transformation (\ref{6: Gauge transf in MaC}) associated to a degree zero
element
$$
\mathsf{H}=\left(H\in \Def(\cA ss_\infty\rar f\cG raphs),\ \ \ h\in  \mathsf{GC}_2\right)
$$
preserves the subset $\cM\cC_{\cA ss}\subset \cM\cC$ if and only if $\delta_\bubu\ga=0$,
i.e.\ if $\ga$ is a cycle in the Kontsevich graph complex. In this case one has
$$
 e^{\mathrm{ad}_{\mathsf{H}}} (\Ga, \bubu) = \left(e^{\mathrm{ad}_{(H\circ e^h)}}(\Ga\circ (e^{-h})) + \ldots\ ,\  \bubu\right)
$$
where $e^{\mathrm{ad}_{(a\circ e^h)}}$ is computed with respect to the Lie bracket in
$\Def(\cA ss_\infty\hspace{-1mm}\rar\hspace{-1mm} \cG raphs)$ and, for an element
$A\in \Def(\cA ss_\infty\hspace{-1mm}\rar\hspace{-1mm} \cG raphs)$ and an element
$\ga\in (\mathsf{GC}_2$ we set
$$
A\circ (e^\ga):=
\sum_{n=0}^\infty \frac{1}{n!} (\ldots ((A\circ \underbrace{\ga)\circ \ga)\ldots \circ\ga)}_n \in
\Def(\cA ss_\infty\hspace{-1mm}\rar\hspace{-1mm} f\cG raphs)
$$
It is clear from  these formulae that gauge transformations of the set $\cM\cC_{\cA ss}$
associated with degree zero elements in $\mathsf{MaC}(\fW)$  of the form
$$
\mathsf{H}=\left(H,0\right)
$$
will induce ---  via the projection(\ref{5: surj of mapping cones})) ---  $\cA ss_\infty$ isomorphisms  of $\cA ss_\infty$ structures associated to elements of $\cM\cC_{\cA ss}$.

\subsection{A naive action of $GRT$ on $\cM\cC_{\cA ss}$} For any $[\ga]\in \fg\fr\ft$ an associated
degree zero element $\mathsf{H}_\ga=(0,\ga)\in \mathsf{MaC}(\fW)$ gives us a gauge transformation of
$\cM\cC$ which preserves the subset $\cM\cC_{\cA ss}$. For example, in the case of the standard Gerstenhaber algebra structure (\ref{4: Ga_0}) one has
\Beq\label{6: g action on Ga_0}
\mathsf{\Ga_0}^\ga:= e^{\mathrm{ad}_{\mathsf{H}_\ga}}\mathsf{\Ga_0}=
 \left( \xy
 (0,0)*{\circ}="a",
(3,0)*{\circ}="b",
\endxy\ +\ \xy  (0,6)*{^g},
 (0,4)*{\bu}="a", (0,0)*{\circ}="b",
\ar @{-} "a";"b" <0pt>
\endxy\ \ , \ \  \xy
 (0,0)*{\bullet}="a",
(4,0)*{\bu}="b",
\ar @{-} "a";"b" <0pt>
\endxy\right)
\Eeq
where $g=\exp(-\ga)\in GRT$. The associated (via the projection (\ref{5: surj of mapping cones})) $\cA ss_\infty$ structure on polyvector fields is precisely the standard differential Gerstenhaber algebra structure in which the differential is twisted by the action of $g$ on the Poisson structure (see Main Theorem in the introduction).

\sip

To construct a less naive action of $GRT$ on  $\cM\cC_{\cA ss}$ we need some technical preparations.

\subsection{Splitting of the Lie algebra  $\sM\sa\sC(\fW)$}

The natural epimorphism of differential Lie algebras,
$$
\sM\sa\sC(\fW) \lon \mathsf{GC}_2,
$$
has a section in the category of {\em non}-differential Lie algebras given explicitly
in the following proposition.

\begin{proposition} There is a morphism of Lie algebras $s: \mathsf{GC}_2 \rar \sM\sa\sC(\fW)$
given by
$$
\Ba{rccc}
\fs:&   \mathsf{GC}_2 & \lon & \Def(\cA ss_\infty\hspace{-1mm}\rar\hspace{-1mm} f\cG raphs)
\oplus \mathsf{GC}_2    \\
& \ga  & \lon & (\ga_\circ, \ga)
\Ea
$$
where
$$
\ga_\circ:=\sum_{v\in V(\ga)}\ga_{v\rar \circ}
$$
and $\ga_{v\rar \circ}$ stands for the graph
$\ga$ whose (black) vertex $v$ is made white .
\end{proposition}\label{6: prop on s}
\begin{proof} Denoting $\var:=|\ga_1||\ga_2|$, we have
\Beqrn
\fs([\ga_1,\ga_2]) & =& \left( \sum_{v\in V([\ga_1,\ga_2])}[\ga_1,\ga_2]_{v\rar \circ}\ \ , \ \
[\ga_1,\ga_2]      \right)\\
&=& \left(\sum_{w\in V(\ga_1)}\left(\sum_{v\in V(\ga_2)} \ga_1\circ_w (\ga_2)_{v\rar \circ} +
\sum_{v\in \{V(\ga_1)\setminus w\}} (\ga_1)_{v\rar \circ}\circ_w \ga_2
\right) -
 (-1)^{\var}(1\leftrightarrow 2)\
,\ [\ga_1,\ga_2]      \right)\\
&=& \left(\sum_{w\in V(\ga_1) \atop v\in V(\ga_2)} \ga_1\circ_w (\ga_2)_{v\rar \circ}
-  (-1)^{\var} \sum_{w\in V(\ga_2) \atop v\in V(\ga_1)} \ga_2\circ_w (\ga_1)_{v\rar \circ}
 + (\ga_1)_\circ \circ \ga_2 -
  (-1)^{\var}(\ga_2)_\circ \circ \ga_1 \
,\ [\ga_1,\ga_2]      \right)\\
&=& \left([(\ga_1)_\circ, (\ga_2)_\circ]  + (\ga_1)_\circ \circ \ga_2 -
  (-1)^{\var}(\ga_2)_\circ \circ \ga_1 \
,\ [\ga_1,\ga_2]      \right)\\
&=& \left[\fs(\ga_1), \fs(\ga_2)\right].
\Eeqrn
\end{proof}

\begin{corollary} There is an isomorphism of Lie algebras
$$
\Ba{rccc}
\fs:&  \sM\sa\sC(\fW) & \lon & \Def(\cA ss_\infty\hspace{-1mm}\rar\hspace{-1mm} \cG raphs)
\oplus \mathsf{GC}_2    \\
& (a,\ga)  & \lon & (a-\ga_\circ, \ga).
\Ea
$$
and hence an isomorphism of gauge groups,
$$
e^{\sM\sa\sC(\fW)^0}\simeq e^{  \Def(\cA ss_\infty\hspace{-1mm}\rar\hspace{-0.3mm} f\cG raphs)^0}\times
e^{ \mathsf{GC}_2^0}.
$$

\end{corollary}

\bip

Consider now an action of $GRT$ on $\Ga_0\in \cM\cC_{\cA ss}$ via the morphism $\fs$,
\Beqrn
e^{\mathrm{ad}_{\fs(\ga)}}\mathsf{\Ga_0}&=&
 \left( \xy
 (0,0)*{\circ}="a",
(3,0)*{\circ}="b",
\endxy + \xy
 (0,4)*{\bu}="a", (0,0)*{\circ}="b",
\ar @{-} "a";"b" <0pt>
\endxy  +  \left[\ga_\circ,\xy
 (0,0)*{\circ}="a",
(3,0)*{\circ}="b",
\endxy + \xy
 (0,4)*{\bu}="a", (0,0)*{\circ}="b",
\ar @{-} "a";"b" <0pt>
\endxy\right] + \ga_\circ \circ \bubu -  \xy
 (0,4)*{\bu}="a", (0,0)*{\circ}="b",
\ar @{-} "a";"b" <0pt>
\endxy\circ \ga + \f(\ga^2)
 \ \ , \ \  \xy
 (0,0)*{\bullet}="a",
(4,0)*{\bu}="b",
\ar @{-} "a";"b" <0pt>
\endxy\right)\\
&=&
 \left(\xy
 (0,0)*{\circ}="a",
(3,0)*{\circ}="b",
\endxy + \xy
 (0,4)*{\bu}="a", (0,0)*{\circ}="b",
\ar @{-} "a";"b" <0pt>
\endxy -\delta_\wiwi\ga_\circ +
\xy (0,6)*{^{\ga}},
 (0,4)*{\bu}="a", (0,0)*{\circ}="b",
\ar @{-} "a";"b" <0pt>\endxy -
\xy (0,6)*{^{\ga}},
 (0,4)*{\bu}="a", (0,0)*{\circ}="b",
\ar @{-} "a";"b" <0pt>\endxy +
\f(\ga^2)  \ \ , \ \  \xy
 (0,0)*{\bullet}="a",
(4,0)*{\bu}="b",
\ar @{-} "a";"b" <0pt>
\endxy\right)\\
&=&
 \left(\xy
 (0,0)*{\circ}="a",
(3,0)*{\circ}="b",
\endxy + \xy
 (0,4)*{\bu}="a", (0,0)*{\circ}="b",
\ar @{-} "a";"b" <0pt>
\endxy -\delta_\wiwi\ga_\circ +
\f(\ga^2)  \ \ , \ \  \xy
 (0,0)*{\bullet}="a",
(4,0)*{\bu}="b",
\ar @{-} "a";"b" <0pt>
\endxy\right)\\
\Eeqrn
As terms of the form $\xy (0,6)*{^{\ga}},
 (0,4)*{\bu}="a", (0,0)*{\circ}="b",
\ar @{-} "a";"b" <0pt>\endxy$ cancel out, the $\cA ss_\infty$ structure on polyvector fields corresponding
to $e^{\mathrm{ad}_{\fs(\ga)}}\mathsf{\Ga_0}$ has the differential,  $\xy
 (0,4)*{\bu}="a", (0,0)*{\circ}="b",
\ar @{-} "a";"b" <0pt>\endxy$, unchanged by the action of $GRT$ at the price of adding higher
homotopies to the standard wedge product. This rather unusual universal $\cA ss_\infty$ structure
is $\cA ss_\infty$ isomorphic to the naive $GRT$ deformation (\ref{6: g action on Ga_0})
since
$$
e^{\mathrm{ad}_{\fs(\ga)}}\mathsf{\Ga_0}= e^{\mathrm{ad}_{\fs(\ga)}} e^{-\mathrm{ad}_{\mathsf{H}_\ga}}\mathsf{\Ga}_0^\ga
$$
and $e^{\mathrm{ad}_{\fs(\ga)}} e^{-\mathrm{ad}_{\mathsf{H}_\ga}}$ is of the form $e^{\mathrm{ad}_{\mathsf{H}}}$ for some
$\mathsf{H}=(H\in \Def(\cA ss_\infty\rar f\cG raphs), 0)$. However this fact does not prove our Main Theorem as the multiplication operation in the $\cA ss_\infty$ algebra corresponding to
$e^{\mathrm{ad}_{\fs(\ga)}}\mathsf{\Ga_0}$ is given by the graph
$$
\xy
 (0,0)*{\circ}="a",
(3,0)*{\circ}="b",
\endxy  - \delta_\wiwi \ga_\circ + \f(\ga^2)
$$
and hence is
{\em not}\, equal to the standard wedge product. However it is now clear how to achieve a $GRT$ deformation of the standard dg algebra structure on polyvector fields in such a way that the differential and the wedge product stay unchanged. In the notations of Lemma~{\ref{5: Lemma on reprent of grt witn no black vertices}}, consider a degree zero map,
given by
$$
\Ba{rccc}
\hat{\fs}:&   \mathsf{GC}_2 & \lon & \Def(\cA ss_\infty\hspace{-1mm}\rar\hspace{-1mm} f\cG raphs)
\oplus \mathsf{GC}_2    \\
& \ga  & \lon & (\ga_\circ+ \ga_\wiwi +\ldots + \ga^{max}_{\circ\ldots\circ}, \ga).
\Ea
$$
Then, for $\ga$ a cycle in $\mathsf{GC}_2$   representing some cohomology class $[\ga]\in \fg\fr\ft$, we have
\Beq\label{6: shat of Ga_0}
e^{ad_{\hat{\fs}(\ga)}}\mathsf{\Ga_0}= \left(\xy
 (0,0)*{\circ}="a",
(3,0)*{\circ}="b",
\endxy + \xy
 (0,4)*{\bu}="a", (0,0)*{\circ}="b",
\ar @{-} "a";"b" <0pt>
\endxy -\delta_\wiwi\ga_\circ +
\f(\ga^2)  \ \ , \ \  \xy
 (0,0)*{\bullet}="a",
(4,0)*{\bu}="b",
\ar @{-} "a";"b" <0pt>
\endxy\right)
\Eeq
so that the first corrections to the standard wedge multiplication, $\wiwi$, in polyvector fields is given by the following graph
$$
\xy
 (0,0)*{\circ}="a",
(3,0)*{\circ}="b",
\endxy  - \delta_\wiwi \ga^{max}_{\circ\ldots\circ} + \f(\ga^2)
$$
As  $\ga^{max}_{\circ\ldots\circ}$ has at least four white vertices, we conclude that
the universal $\cA ss_\infty$ structure corresponding to $e^{\mathrm{ad}_{\hat{\fs}(\ga)}}\mathsf{\Ga_0}$ has operations $\mu_1$ and $\mu_2$ unchanged at the price of non-trivial higher homotopy operations $\mu_{\bu \geq 4}\neq 0$. We have
$$
e^{\mathrm{ad}_{\hat{\fs}(\ga)}}\mathsf{\Ga_0}= e^{\mathrm{ad}_{\hat{\fs}(\ga)}} e^{-\mathrm{ad}_{\mathsf{H}_\ga}}\mathsf{\Ga}_0^\ga= e^{-\mathrm{ad}_{(H,0)}}\mathsf{\Ga}_0^\ga
$$
for some $H\in \Def(\cA ss_\infty\rar f\cG raphs)$. Thus the universal $\cA ss_\infty$ structures corresponding to Maurer-Cartan elements (\ref{6: g action on Ga_0}) and (\ref{6: shat of Ga_0})
are $\cA ss_\infty$ isomorphic. This proves our Main Theorem for the case $M=\R^d$, the affine space.
\subsubsection{{\bf Globalization to any Poisson manifold}}
Let \(M\) be a finite-dimensional smooth manifold. A torsion-free affine connection on \(M\) defines an isomorphism of sheaves of algebras between the sheaf of jets of functions, \(J^{\infty}C^{\infty}_M\), and the completed symmetric bundle \(\hat{S}(T^*_M)\) of the cotangent bundle. Similarly, the sheaf of jets of polyvector fields, \(J^{\infty}(S(T_M[-1]))\), becomes isomorphic to the sheaf \(\mathfrak{T}:=\hat{S}(T^*_M\oplus T_M[-1])\). The canonical jet bundle connection defines, \textit{via} this isomorphism, a Maurer-Cartan element \(B\in\Omega(M,\mathfrak{T})\) of the dg Lie algebra of differential forms on \(M\) with values in \(\mathfrak{T}\). Taking jets (with respect to the affine connection) is a quasi-isomorphism
 \[
  j:(\cT_{poly}(M),\wedge,[\,,\,]_S)\hookrightarrow (\Omega(M,\mathfrak{T}),d_{dR}+[B,\,]_S,\wedge,[\,,\,]_S)
 \]
of Gerstenhaber algebras. The space on the right was used, e.g., in \cite{D}, to globalize Kontsevich's formality morphism. The action of degree \(0\) cocycles of Kontsevich's graph complex \(\mathfrak{GC}_2\) by \(\mathcal{L}ie_{\infty}\)-derivations of the polyvector fields on affine \(\R^d\) defines, because of equivariance with respect to linear coordinate changes, \(\mathcal{L}ie_{\infty}\)-derivations of the dg Lie algebra \((\Omega(M,\mathfrak{T}),d_{dR},[\,,\,]_S)\). Let now \(\pi\) be a Poisson bivector on \(M\). The jet \(j(\pi)\) is then a Maurer-Cartan element of \(\Omega(M,\mathfrak{T})\) and, because \(\mathcal{L}ie_{\infty}\) morphisms can be twisted by Maurer-Cartan elements, any degree \(0\) graph cocycle \(\gamma\) will define a \(\mathcal{L}ie_{\infty}\) morphism
 \[
  e^{\gamma}:(\Omega(M,\mathfrak{T}),d_{dR}+[j(\pi)+B,\,]_S,[\,,\,]_S)\rightarrow
(\Omega(M,\mathfrak{T}),d_{dR}+[\gamma(j(\pi)+B),\,]_S,[\,,\,]_S).
 \]
Define \(\delta:=d_{dR}+[B,\,]_S\). Because \(\gamma(j(\pi)+B)=j(\gamma(\pi))+B\), the \(\gamma\) on the right referring to the globalized automorphism of polyvector fields on \(M\) (see \cite{D} for the arguments), the above is a morphism
 \[
  (\Omega(M,\mathfrak{T}),\delta+d_{j(\pi)},[\,,\,]_S)\rightarrow(\Omega(M,\mathfrak{T}),\delta+d_{j(\gamma(\pi))},[\,,\,]_S).
 \]
Our formula for the morphism of associative Poisson cohomology algebras on affine space is invariant under linear coordinate changes so it defines a morphism of associative algebras
 \[
  F^{\gamma}:H((\Omega(M,\mathfrak{T}),\delta,\wedge),d_{j(\pi)})\rightarrow H((\Omega(M,\mathfrak{T}),\delta,\wedge),d_{j(\gamma(\pi))}).
 \]
Since taking jets is a quasi-isomorphism of associative algebras,
 \[
  (\cT_{poly}(M),d_{\pi},\wedge)\rightarrow (\Omega(M,\mathfrak{T}),\delta+d_{j(\pi)},\wedge,),
 \]
this shows that the morphism \(F^{\gamma}\) globalizes.

\bip

\bip

\bip

\def\cprime{$'$}

\end{document}